\documentclass{article}
\usepackage{graphicx}
\usepackage{epstopdf}

\usepackage{hyperref}
\usepackage{comment}

\title{Symmetrisation of $n$-operads and compactification of real configuration spaces}
\author{M.A. Batanin\protect \footnote{The author holds the Scott Russell Johnson Fellowship in
the Centre of Australian Category Theory at Macquarie University}\\ Macquarie University,  NSW 2109, Australia \\
e-mail: \href{mailto:mbatanin@math.mq.edu.au}{mbatanin@math.mq.edu.au}} 

\date{March 15, 2006}

\newtheorem{theorem}{\bf Theorem}[section]
\newtheorem{defin}{\bf Definition}[section]
\newtheorem{pro}{\bf Proposition}[section]
\newtheorem{lemma}{\bf Lemma}[section]
\newtheorem{corol}{\bf Corollary}[pro]
\newtheorem{cor}{\bf Corollary}[theorem]

\newcommand{\F}{\mbox{$\cal F$}}
\newcommand{\PF}{\mbox{$\cal PF$}}
\newcommand{\RF}{\mbox{$\cal RF$}}
\newcommand{\URF}{\mbox{${\cal RF}_n^{\circ}$}}
\newcommand{\RSF}{\mbox{${\cal RF}_{\infty}$}}

\newcommand{\Q}{   \hfill  \ $\scriptstyle \clubsuit$ \ }
\newcommand{\h}{\mbox{$\bf h$}}

\newcommand{\PH}{\mbox{$\bf PH$}}
\newcommand{\RH}{\mbox{$\bf RH$}}
\newcommand{\ph}{\mbox{$\bf ph$}}

\newcommand{\rh}{\mbox{$\bf rh$}}
\newcommand{\scrh}{\mbox{$\bf scrh$}}
\newcommand{\urh}{\mbox{$\mathbf{rh}_{\circ}^n$}}
\newcommand{\URH}{\mbox{$\mathbf{RH}_{\circ}^n$}}
\newcommand{\dom}{\triangleleft}

\renewcommand{\H}{\mbox{$\bf H$}}
\newcommand{\fm}{\mbox{$\bf fm$}}
\newcommand{\scfm}{\mbox{$\bf scfm$}}

\newcommand{\GJ}{\mbox{$\bf GJ$}}
\newcommand{\SCGJ}{\mbox{$\bf SCGJ$}}
\newcommand{\J}{\mbox{$\bf J$}}

\newcommand{\M}{\mbox{$\bf m$}}

\newcommand{\FN}[1]{\mbox{$\mathrm{Mod}_{{#1}}^{\ \! n}$}}
\newcommand{\mod}[1]{\mbox{$\mathrm{mod}_{[{#1}]}^{\ \! n}$}}
\newcommand{\Conf}{\mbox{ Conf}}

\newcommand{\fn}{\mbox{$\mathrm{FN}$}}
\newcommand{\PO}{\mbox{${PO}_n^{(n-1)}$}}
\newcommand{\RO}{\mbox{$RO_n^{(n-1)}$}}

\newcommand{\URO}{\mbox{$URO_n^{(n-1)}$}}

\newcommand{\Coll}{\mbox{$ Coll_n^{S}$}}
\newcommand{\Oper}{\mbox{$ Oper_n^{S}$}}

\newcommand{\colim}{\mbox{$\mbox{\rm co}\!\!\lim\limits_{\ph^{n}_k}$}}
\newcommand{\hocolim}{\mbox{$\mbox{\rm hoco}\!\!\lim\limits_{\ph^{n}_k}$}}

\newcommand{\Proof}{\noindent {\bf Proof. \ }}
\renewcommand{\theequation}{\thesection.\arabic{equation}}
\makeatletter\@addtoreset{equation}{section}\makeatother

\newcommand{\Example}{\noindent \makebox[25mm]{{\bf
 Example  \hspace{-3mm}
\addtocounter{example}{1} 
\thesection.
 \hspace{-3.5mm}
\theexample \ \  }}}

\newtheorem{conj}{\bf Conjecture}[section]

\newcommand{\Remark}{\noindent \makebox[23mm]{{\bf Remark
\hspace{-1mm}\addtocounter{remark}{1} 
\thesection.\theremark \ }}}
\newcounter{remark}[section]
\newcounter{example}[section]

\begin{document}
\maketitle

\begin{abstract}

It is well known that the forgetful functor from symmetric operads to nonsymmetric operads has a left adjoint $Sym_1$ given by 
product with the symmetric group operad.   It is also well known that this functor does not affect the category of algebras of the operad. From the point of view of the author's theory of higher operads, the nonsymmmetric operads are $1$-operads and $Sym_1$ is the first term of the infinite series of left adjoint functors $Sym_n,$ called symmetrisation functors, from $n$-operads to symmetric operads with the property that the  category of one object, one arrow , . . . , one $(n-1)$-arrow algebras of an $n$-operad $A$ is isomorphic  to the
category of  algebras
 of $Sym_n(A)$.  
 
 In this paper we consider some geometrical and homotopical aspects of the  symmetrisation of $n$-operads. 
 We follow Getzler and Jones and consider their decomposition  of the Fulton-Macpherson operad of compactified  real configuration spaces. We  construct 
  an $n$-operadic counterpart of this  compactification which we call the Getzler-Jones operad. We study the properties of Getzler-Jones operad and find that it is contractible and cofibrant in an appropriate model category. The symmetrisation of the Getzler-Jones operad turns out to be exactly the operad of Fulton and Macpherson.  These results should be considered as an
extension of Stasheff's  theory of $1$-fold loop spaces to $n$-fold loop spaces $n\ge 2 .$ We also show that a space $X$ with an action of a contractible $n$-operad has a natural structure of an algebra over an operad weakly equivalent to the little $n$-disks operad. A similar result  holds for chain operads.  These results  generalise the classical   Eckman-Hilton argument to arbitrary dimension.

Finally, we apply the techniques to the Swiss Cheese type operads introduced by Voronov and prove analogous results in this case.    

\

1991 Math. Subj. Class.  18D05 , 18D50, 55P48     
\end{abstract}

\pagebreak

\tableofcontents
        
\section{Introduction}

This is the second  paper in a sequence of  papers devoted to the
relations between higher categories and $n$-fold loop space theory. In the
first paper \cite{BEH}
we developed the necessary categorical techniques which allow us to go back
and forth between   $n$-operads and  classical symmetric operads. The main
goal of this paper
is to clarify the geometric and homotopy theoretic aspects of this theory.

 To do this we restrict ourselves  to a class of so called pruned $(n-1)-$terminal
$n-$operads. This is a slightly smaller category  of $n$-operads than we
considered in \cite{BEH} but it is big enough to include most applications
we have in mind. The reason is that the functor of desymmetrisation
from symmetric operads to $n$-operads \cite{BEH} can be factorised through
the category of pruned $n$-operads. Moreover, as in the unpruned case  this desymmetrisation functor  preserves endomorphism operads. This allows us to construct
a theory of symmetrisation very much in parallel to the unpruned case.
It turns out that pruned $n$-operads are easier to handle  from the
combinatorial point of view.
We have a conjecture,
however, that the main results of this paper are also true in the  case of
$(n-1)$-terminal $n$-operads but so far we have been unable to prove it in
this generality. 

We apply the categorical methods of \cite{BEH} to the category of pruned
$(n-1)$-terminal $n$-operads and to the even smaller class of reduced 
$(n-1)$-terminal $n$-operads. These methods immediately imply the existence
of some categorical symmetric operads $\ph^n$ and $\rh^n ,$ which represent
the theories of  internal  $(n-1)$-terminal pruned $n$-operads  and
internal $(n-1)$-terminal reduced $n$-operads inside categorical symmetric
operads, in full analogy with the categorical operad $\h^n$ in \cite{BEH}.
Our first significant result here is   Theorem \ref{phnrhnmn} which asserts that the simplicial 
operads  $N(\ph^n)$ and $N(\rh^n) ,$ where $N$  is the  nerve functor, are   $E_n$-operads in the
category of simplicial  symmetric operads and reduced simplicial symmetric
operads respectively.

Together with the theorems \ref{formulap}  and  \ref{formular}, these theorems show
that once we have a space with an action of a contractible pruned
$(n-1)$-terminal $n$-operad there is an action of an $E_n$-operad on this
space. An analogous result holds for the reduced operads in chain complexes.  As we conjectured in \cite{BEH} this
should give a very natural proof of Deligne's conjecture answering a question by  
Kontsevich \cite{K}\footnote{During the preparation of this paper,
 D.Tamarkin informed me that he indeed obtained a proof of Deligne's conjecture
by exhibiting a contractible $2$-operad acting naturally on a  $2$-graph consisting of  $DG$-categories, $DG$-functors and the  complex of their derived natural transformations. In the particular case of a $DG$-category with one object and transformations  of the identity functor, we obtain the  Hochschild complex of an  associative algebra and then we  apply  our Theorem \ref{Enchain}
 \cite{Tam1}.}.

There are, however, other results in the present paper which we believe are significant.
First we observe that the desymetrisation of the operad of Fulton and
Macpherson $\fm^n$ of compactified real configuration spaces  
\cite{GJ,KS} contains a contractible reduced $(n-1)$-terminal $n$-operad 
which we call the Getzler-Jones operad $\GJ^n .$  Actually, this operad was
discovered by Getzler and Jones in their remarkable preprint \cite{GJ}.
The apparatus of $n$-operads did not exist at that time and Getzler and
Jones attempted to express the properties of $\GJ^n$ in terms of a natural 
subdivision of the operad $\fm^n$. It turned out, however, that   Getzler-Jones subdivision does not
give a cellular structure compatible with the operadic structure of $\fm^n $ as was first observed by
D.Tamarkin (see \cite{VoronovG} for an explanation of Tamarkin's
counterexample.) This counterexample implied considerable technical
difficulties in the proof of Deligne's conjecture.

The second implicit appearance of $\GJ^n$ was in the Kontsevich and Soibleman
paper \cite{KS}. For $n=2$ they considered  closed
contractible subsets $X_T$ of $\fm^2$  where $T$ is a finite set with two
complementary orders on it (see Definition \ref{complorder}). They used some
properties of $X_T$ to prove  Deligne's conjecture for
$A_{\infty}$-algebras.
We show that the  space $\GJ^n_T$  of arity $T$  is indeed  equal to the
Kontsevich-Soibelman $X_T$ for a pruned $n$-tree $T$ (see Section
\ref{Comporders} for the explanation  of the  connection between trees
and complementary orders, and Section \ref{GJoperad} for the definition
of $X_T$ for any $n$).

We show that $\GJ^n$ can be considered as a natural analogue of $\fm^n$ in
the category of  reduced $(n-1)$-terminal $n$-operads. In particular, set
theoretically it is a free reduced $(n-1)$-terminal $n$-operad on the reduced 
$n$-collection of Getzler-Jones cells. This last result leads to the
theorem that the symmetrisation of $\GJ^n$ is exactly $\fm^n$ which is the
basis for the main results of our paper. We also think that this gives an
interesting new insight to the geometry of $\fm^n .$ We are going to
continue this study in the next paper of this series.

The operad $\GJ^n$ is still not cellular in the strong sense that 
 it is not a geometric realisation of a poset of cells of a regular CW-complex. 
 Tamarkin's counterexample works
well in this case too.  Nevertheless, its unbased version $\GJ^n_{\circ}$ is  a cellular object in the category
of unbased reduced $(n-1)$-terminal $n$-operads in the model category theoretic sense, and
in particular it is a cofibrant contractible operad. The term  unbased here means  that we forget
 about nullary operations of our operads. There is also an  $n$-operad  $\URH,$ which is an unbased 
categorical reduced $n$-operad freely generated by its internal reduced $n$-operad. The geometric 
realisation of the  nerve of this operad is  cofibrant and contractible and  is strongly
homotopy equivalent to the operad $\GJ^n_{\circ} .$  This implies a homotopy equivalence between the geometric realisation of 
 $N(\urh)$ and $\fm^n_{\circ}.$ This  operad $\urh$ (the unbased 
categorical reduced  symmetric operad freely generated by its internal reduced $n$-operad)
  plays, therefore,  the role of the
nonexistent poset operad of cells of the Getzler and Jones decomposition.

Perhaps, the most intresting result of the above study of the combinatorics of the Getzler-Jones operad 
is the following generalisation of Stasheff's classical theory of $A_{\infty}$-spaces.
The closure of a Getzler-Jones cell $K_T= cl(\FN{T})$ inside $\GJ_T$ is a 
manifold with corners homeomorphic  to the ball of dimension $E(T)-n-1 ,$ where $E(T)$ is the number 
of edges of the $n$-tree $T .$ 
The cellularity   of $\GJ^n_{\circ}$ means that an action of $\GJ^n_{\circ}$ 
on a pointed space $X$ can be described as an inductive process of extension 
of  higher   homotopies from the boundary of $K_T$ to its interior in exact 
analogy with Stasheff's  description of  $A_{\infty}$-spaces. 
 And our Theorem \ref{SymGJ} states that  the collection $K_{\bullet},$ 
where $\bullet$ runs over the set of pruned $n$-trees,  gives  full coherence 
conditions for $E_n$-spaces. If $n=1$  the collection  $K_{\bullet}$ is 
the sequence of associahedra \cite{GJ} and  we get  Stasheff's theorem. 

The difference between $1$-dimensional  and higher dimensional  cases 
appears in the existence of some cells (Tamarkin's cells) in $\GJ^n$
 which are not completely on  the boundary of $K_T.$ This is not  excessive
 information, however,  just a defect of our language when we try to express 
the coherence laws in terms of an action of an $(n-1)$-terminal $n$-operad. 
There is a way of avoiding  this problem  but we will have to pay a price by
 using more sophisticated $n$-operads which have full source and target operations. 
We will consider this subject in a forthcoming paper.

Finally, in the last section of this paper we apply our techniques to the
case of Swiss Cheese type operads \cite{VoronovSC}. The advantage of our
categorical  methods is that we have nothing to prove here once we put
the right definitions of our main categories and functors in place. We deduce
immediately a symmetrisation formula for Swiss Cheese type $n$-operads
\footnote{I am  grateful to D.Tamarkin for encouraging me  to look at the action of the Swiss
 Cheese operad from the $n$-operadic 
point of view and for sending me a preliminary version of his papers concerning the action 
of Swiss-Cheese operads on associative algebras and their  Hochschild complexes \cite{Tam2}.}
  and  other Swiss-Cheese analogues of the results for classical operads.

We hope  that similar results can be obtained for symmetrisation of 
some other important coloured operads; for example, operads for morphisms
between $E_n$-algebras. This should lead to a better understanding of
coherence conditions for   such morphisms.

\

\Remark  We will freely use  the terminology from \cite{BEH}  concerning $n$-trees,
 $(n-1)$-terminal $n$-operads and their
algebras. Notice, however, that the main objects for us here are $(n-1)$-terminal 
$n$-algebras of our $n$-operads. Roughly speaking such an algebra is an object $X$ 
of our basic symmetric monoidal category $V$ together with an action 
$A_T\otimes X^{\otimes^k} \rightarrow X$ where $k$ is the number of tips of 
$T$ satisfying some natural conditions.  An $n$-operad $A$ may have, however,
 more complicated types of algebras which involve source and target operations but we do not need them here. 
So we will speak simply about category of algebras of $A$ having in mind the 
subcategory of its $(n-1)$-terminal algebras. We refer the reader to \cite{BEH} for more discussion about this issue. 

\

\noindent {\bf Acknowledgements.}  I would like to thank Andr\'{e} Joyal  for pointing 
out a mistake in the first version of this paper. It is exactly the  work on this mistake 
 which led me to an understanding of the role of the Getzler-Jones operad in $n$-operad  theory. 
I am also very grateful to him for organising my visit to Montreal in July 2003 and for our 
fruitful and inspiring discussions during that time.  

 I  also wish to express my  gratitude to
C.Berger, L.Breen,  E.Getzler,   S.Lack,  M.Markl,  
  J.Stasheff,  R.Street,  D.Tamarkin and A.Voronov  
for numerous illuminating conversations and correspondences. 
 
 Finally, I gratefully acknowledge the financial support of the Scott Russell Johnson Memorial 
 Foundation and Australian Research Council (grant \# DP0558372).

 \section{Complementary orders and  pruned $n$-trees}\label{Comporders}
 
 Here we consider a combinatorial techniques which will be used further to develop a theory of pruned $n$-operads. 
The machinery of bar-codes from  \cite{GJ}  is an equivalent language but  we prefer to work with  the
 notion of complementary orders   introduced by Kontsevich and Soibelman in \cite{KS}. 
 
 \begin{defin} A partial  order on a set $X$ is a transitive, 
antireflective relation on $X .$ It is called linear if any two elements are comparable.
 \end{defin}  
 Antireflective here means  that the diagonal is always in the complement relation.
 We choose this terminology just for  technical reasons. 
Of course, we always can add the diagonal to our relations 
and work with reflective relations.  So if we are given 
a partial order $<$ we will often use the notation $a\le b$ to mean that 
$a<b $ or $a=b .$ This does not lead to any trouble.  
   Observe, also that  a partial order in our sense is always antisymmetric.

 \begin{defin}[Kontsevich-Soibelman \cite{KS}]\label{complorder} 
 Let $I$ be a  set. Suppose we have an $n$-tuple of  
 partial orders $\Xi = (<_0, \ldots, <_{n-1})$ on $I .$ 
We call them complementary  orders provided any two elements 
 $i, j\in I$ can be compared with respect to exactly one order  $<_0, \ldots, <_{n-1} .$ \end{defin}
 \begin{defin} A set with a given  $n$-tuple of  complementary orders on it will be called an $n$-ordered set. \end{defin}

\begin{lemma} Let  $T$ be a pruned $n$-tree
$$T=
[k_n]\stackrel{\rho_{n-1}}{\longrightarrow}[k_{n-1}]\stackrel{\rho_{n-2}}{\longrightarrow}...
\stackrel{\rho_0}{\longrightarrow} [1] .
$$  
Then the relation: $i <_{p} j$  
 if and only if   
 $i<j$
 in $[k_n]$ and
 $$ \rho_{n-1}\cdot ... \cdot \rho_{p}(i) = \rho_{n-1}\cdot
\ldots
\cdot \rho_{p}(j)$$ but 
$$ \rho_{n-1}\cdot ... \cdot \rho_{p+1}(i) \ne \rho_{n-1}\cdot \ldots
\cdot \rho_{p+1}(j) ;$$
defines $n$ complementary orders on the set of tips of $T .$
\end{lemma}

\Proof The proof is obvious. 

\

\Q 

\

In fact we can characterise pruned $n$-trees in the following way:

\begin{lemma}\label{treeord} Suppose we are given $n$-complementary orders $<_0, \ldots, <_{n-1}$ on  
 a finite set $X$  such that 
\begin{itemize}
\item  if $i<_p j$ and $j<_r l$ then $i<_{min(r,p)} j$
\end{itemize}
then there exist a linear order on $X$ and  a pruned $n$-tree $T$ such that the ordinal of 
its tips is $X \simeq [k_n]$ and the complementary orders $<_0, \ldots, <_{n-1}$ are determined by $T .$
\end{lemma}

\Proof  The linear order $<$ on $X$ is defined by the requirement that $i<j$ if and only 
if there exists a $p$ such that $i<_p j .$ Let $[k_n]$ be the corresponding ordinal. 

Suppose that there exists a triple $i<j<l$ from $[k_n]$ such that  $i<_{n-1}<l$ but $i<_r j$ for
 some $r< n-1$ then $i<_r j$ so $i<_r l ;$ contradiction. Similarly for the other side.
 So the order $<_{n-1}$ determines a subdivision on $k_{n-1}$ intervals of the ordinal $[k_n] .$ 
 This can be considered as a surjection of ordinals $[k_n]\rightarrow [k_{n-1}] .$ 
Obviously we have $(n-1)$ complementary orders on $[k_{n-1}]$ which satisfy the conditions of the lemma.
 So we can proceed and construct a tree $T .$

\Q 

\begin{defin} We will call an $n$-ordered   set $X$ totally $n$-ordered if
 it satisfies the conditions of lemma \ref{treeord}.  
 If $X = \{1,\ldots,k\}$ is totally $n$-ordered and the induced linear 
order makes it equal to the ordinal $[k]$ then we call $X$  an $n$-ordinal.  
This includes the case of a terminal ordinal (one element set and empty $n$ complementary  orders) and 
of an initial ordinal (empty set). \end{defin}   

\

\Remark We will often consider a special $n$-ordinal for  which only one order $<_l$ is nonempty.
 We will use the notation $M_l^k$ for such an ordinal on $\{1,\ldots,k\}$ (see \cite{BEH} 
for the picture of corresponding pruned tree with the same notation).

\begin{defin} Let $X$ and $Y$ be two $n$-ordered sets. An order preserving map (or  a map of $n$-ordered sets) 
 from $X$ to $Y$ is a map $f:X\rightarrow Y$ such that  $i<_p j$ in $X$   implies that
$f(i) \le_r f(j)$ for some $r\ge p$ or $f(j)<_r f(i)$ for $r>p .$ \end{defin}

\begin{defin} Suppose we are given two $n$-tuples of  complementary orders $\Xi_1$ and $\Xi_2$ 
on the same  set $X.$ We will say that $\Xi_1$ dominates $\Xi_2$ (notation $\Xi_2 \dom \Xi_1$)
 if   $i <_p j$ in $\Xi_1$ implies  $i <_r j$ for some $r\ge p$ or $j<_r i$ for $r>p $ in $\Xi_2 .$ \end{defin} 
Of course, $\Xi_1$ dominates $\Xi_2$ if and only if  the identity map $X \rightarrow X$ is a map of $n$-ordered sets. 

\begin{lemma} Let 
$$T=
[k_n]\stackrel{\rho_{n-1}}{\longrightarrow}[k_{n-1}]\stackrel{\rho_{n-2}}{\longrightarrow}...
\stackrel{\rho_0}{\longrightarrow} [1]
$$  
 and 
 $$S=
[s_n]\stackrel{\xi_{n-1}}{\longrightarrow}[s_{n-1}]\stackrel{\xi_{n-2}}{\longrightarrow}...
\stackrel{\xi_0}{\longrightarrow} [1]
$$  
  be two pruned $n$-trees. Let $f: [k_n] \rightarrow [s_n]$ be a map. Then there exists a map  of trees $\sigma: T \rightarrow S$ such that 
$\sigma_n = f$ if and only if $f$ is a map of $n$-ordered sets.  \end{lemma}

\Proof      If such a $\sigma$ exists the orderpreserving property of $f$ is obvious. 
Now we want to reconstruct a $\sigma$ from $f$. 
 We put $\sigma_n = f: [k_n] \rightarrow [s_n] ,$  of course. 
 Now take a point $a$ from $[k_{n-1}]$ then its  preimage under  
$\rho_{n-1}$ consists of an interval and we define $\sigma_{n-1}(a) = \xi_{n-1}(\sigma_n(i)), $ 
 where $i$ is an arbitrary element from the preimage 
 of $a . $ Since $f$ preserves  order, this definition is correct. 
Indeed, if $i<_{n-1} j$ is another element from the preimage then 
$f(i)\le_{n-1} f(j)$ and, hence, their images under $\xi_{n-1}$ are equal. 
 Now we can check quite easily that the constructed $\sigma_{n-1}$
 preserves the $(n-1)$ complementary orders determined by $\partial S$ 
and $\partial T$ and we can proceed with our construction. 
\nopagebreak[4] 

\

 \Q   
 
 \

 We assume that the only degenerate pruned $n$-tree is  $z^n U_0$ .
We thus have   
\begin{theorem}\label{ordn} The category $Ord_n$ of $n$-ordinals and their
order preserving maps
is isomorphic to the category of pruned trees and their morphisms.
\end{theorem}

 One can consider the poset $\J^n_X$ of all total complementary  $n$-orders
on a fixed  set
$X $ with respect to the domination relation.  If $X = \{1, \ldots, k\}$ we
will denote this poset by $\J^n_k .$ The symmetric group $\Sigma_k$ acts naturally on  $\J^n_k .$

  Let $\Upsilon_n(k)$ be the subcategory of $\Omega_n$ whose objects are
pruned $n$-trees with $k$ tips and  whose morphisms are morphisms of
trees which are bijections on tips. We  call such a morphism {\it a
quasibijection}.  Theorem \ref{ordn} implies
\begin{cor}\label{quasibijection} There is a natural isomorphism of
categories
$$\J^n_k/\Sigma_k \rightarrow \Upsilon_n(k) .$$
\end{cor}

\

\Remark The set $\J^n_k$ appeared many times in the literature
\cite{Berger, GJ, BFSV, FN}.
It is isomorphic to the poset of cells  of the classical Fox-Neuwirth
stratification of 
configuration space $\Conf_k(\Re^n)$ (see
  section \ref{FMandGJ}); it  is also isomorphic to the Getzler and Jones
poset of bar-codes
 and to the  {\it Milgram poset} of \cite{BFSV}. We use the last name in
this paper. 

The category $\Upsilon_n(k)$ is isomorphic to   Berger's {\it shuffle category} \cite{Berger}, where it is also shown how to reconstruct $\J^n_k$ as a Grothendieck construction of a {\it shuffle functor}  $\Upsilon_n(k)\rightarrow Cat .$ 

\

\begin{defin} Let $f: T \rightarrow S$ be an order preserving map of $n$-ordinals. 
 For an element $i\in S ,$ the preimage $f^{-1}(i)$ with its natural structure 
of an $n$-ordinal induced from $T$  will be called the fiber of $f$ over $i$. \end{defin}  

Following \cite{GK} and  \cite{BergerM} we define {\it a tree} to be  an isomorphism class  
 of   finite  connected acyclic graphs with a marked 
vertex $v_0$ called the root.  A  vertex which is not the root vertex and 
having valency more than $2$ is called {\it an internal vertex}.  
 The edges of this graph with one open end  are called {\it leaves.} 
Every edge of a tree has, therefore, a target vertex and a source vertex 
provided this edge is not a leaf. Every vertex also has  a set 
of incoming edges and one outcoming edge if the vertex is not a root.  
A monotone path in a tree is a sequence of edges such that 
the target of each edge in the sequence is equal to the source of the next edge. 
The length of this path is the number of edges in the sequence. 
The monotone paths are ordered by inclusion.  Notice that there can be 
several maximal paths in the tree with respect to this order. 
 The {\it length of a  tree} is the maximal length of the maximal monotone  paths in the  tree.

If we choose an incoming edge $e$ at a vertex $v$ we can construct a subtree by considering all vertices  
 and edges which can be connected to $v$ by a monotone path. We call this subtree  {\it a branch} corresponding to $e .$ 

A tree is called {\it labelled} by the set $\{1, \ldots, k \}$ if there is given a bijection from the set of
 the leaves of the tree to $\{1, \ldots, k \} .$ 
 \begin{defin} A  labelled  $n$-planar tree   is a labelled   tree such that for every internal vertex
 $v$ the set of incoming edges is a totally $n$-ordered set. \end{defin}   
 
 Since every totally $n$-ordered set has a canonical linear order   then every  $n$-planar tree has 
   a  canonical structure of a planar tree. So we often will   speak about $n$-planar trees as planar 
trees decorated by $n$-ordinals.

 Every such  tree determines a structure 
 of $n$-ordered set on the set of its labels.  For a vertex $v$ of a planar tree $\tau$ and a 
label $i\in\{1,\ldots,k\},$ let us define 
 $\#_v(i)$ to be the last incoming edge (if it  exists) in the monotone path which connects 
the outcoming edge with label $i$ and $v. $ Then we put $i<_p j$ if $\#_v(i)<_p \#_v (j)$ 
for a vertex $v\in \tau $ , which always exists and is  unique. 
 The domination relation, therefore, induces a partial order on the set of $n$-planar trees 
 labelled by the set $\{1,\ldots,k\} .$ We will use the same notation $\dom$ for this relation.

\section{Pruned $(n-1)$-terminal $n$-operads}\label{prunedoperads} 

Now we can easily give the definition of an $n$-operad based on $n$-ordinals and their morphisms.
 We will, however, show  that this is just a subcategory of the category of $n$-operads which we
 will call the category 
of pruned $(n-1)$-terminal $n$-operads $\PO(V) ,$  where $V$ stands  
for a symmmetric monoidal category over which we consider the operads.  

 We will  call a morphism of trees  $\sigma: T \rightarrow S$ a {\it full injection} if it is  injective,  and  bijective on tips.
 Obviously every fiber  of $\sigma$
is equal to $U_n .$  If $A$ is an $(n-1)$-terminal $n$-operad in $V$ and $\sigma$ is a full
 injection then we have the following composite morphism
$$ A_S \stackrel{\simeq}{\longrightarrow} A_S\otimes I\otimes \ldots \otimes 
I \rightarrow A_S \otimes A_{U_n}\otimes \ldots \otimes A_{U_n} \rightarrow A_T .$$ 

\begin{defin}We call an $n$-operad {\it pruned} if this morphism is an identity for every full injection $\sigma .$ \end{defin}

\Remark In the case of the empty set of fibers of $\sigma$ (i.e. $S = zS'$), we require the morphism $$ A_S 
\rightarrow A_S\otimes I \otimes I \ldots \otimes I \rightarrow A_T$$
to be an identity. This means, that for any  degenerate tree $S$ there is an identification 
$A_S = A_{z^n U_0} .$  Recall that  the only degenerate pruned tree is  $z^n U_0 .$ 

\begin{defin} A pruned ($(n-1)$-terminal) $n$-collection in $V$ is a family of objects $A_T\in V$,
 where $T$ runs over the set of pruned $n$-trees. They form a category   
 $PColl_n(V)$ with respect to termwise morphisms of collections. \end{defin} 

We will denote by $T^{(p)}$ the maximal pruned subtree of a tree $T .$  
Then we can  reformulate the definition of the pruned $n$-operad in the following way:

\begin{lemma} A pruned  $(n-1)$-terminal $n$-operad is given by
a pruned $(n-1)$-terminal $n$-collection $A$ equipped with : \begin{itemize} 
\item  a morphism $e:I \rightarrow A_{U_n} ;$
\item
 a morphism
$$m_{\sigma}: A_S\otimes A_{T_1^{(p)}}\otimes \ldots \otimes A_{T_k^{(p)}} \rightarrow A_T $$
for every morphism of trees $\sigma: T\rightarrow S$ between pruned trees in $\Omega_n .$
\end{itemize}

They must satisfy the usual identities. 
\end{lemma}

This makes it obvious that the category of pruned operads is isomorphic to 
the category of $n$-operads based on $n$-ordinals.  

         \begin{pro}\label{pon}  If $V$ is a cocomplete symmetric monoidal category then the 
forgetful functor $$PU_n:\PO(V) \rightarrow PColl_n(V)$$ has a  left adjoint $\PF_n$  and is monadic. 

          The free pruned $n$-operad monad on the category of pruned $n$-collections in $Set$  is finitary and cartesian. \end{pro}
 
 \
 
 \Remark By slightly abusing notation we will denote the free pruned $n$-operad monad as well as its functor part by    $\PF_n .$ 
         
         \
         
  \Proof The description of the required monad   on a pruned collection $X$ is analogous to the description of the 
free $n$-operad functor.   It is given by an obvious inductive process. 
  
 Let us call a  pruned tree $T$  an {\it admissible} expression of arity $T .$  We also have an admissible expression 
$e$ of arity $U_n .$   If $\sigma:T\rightarrow S$ is a morphism of pruned trees and  the admissible expressions  
 $x  , x_{1},\ldots,  x_{k}$ of arities $S, T_1^{(p)}, \ldots, T_k^{(p)}$ respectively  are already constructed 
then the expression 
$\mu_{\sigma}(x; x_{1}, \ldots x_{k})$ is also an admissible expression of arity $T .$
We also introduce an obvious equivalence relation on the set of admissible expressions generated by pairs 
of composable morphisms of pruned trees  and by two equivalences $T \sim \mu(T;e, \ldots, e) \sim \mu(e;T) $ 
generated by the identity morphism of $T$ and the unique morphism $T \rightarrow U_n .$ Notice however, 
that there are morphisms of trees all of whose fibers (after the  pruning operation) are equal to $U_n .$ 
 We can form an admissible expression 
$\mu_{\sigma}(S;e, \ldots, e)$  corresponding to such a morphism but it is not equivalent to $S ,$ unless  
$\sigma$ is equal to the identity.

Everything else is in complete analogy with the case of free $(n-1)$-terminal $n$-operads \cite{BEH}.

\
  
  \Q
  
  \
  
   Applying the general theory of internal algebras \cite{BEH} we get 
 
 \begin{corol}\label{PH}       
The $2$-functor  of internal pruned $n$-operads is representable by a pruned  $n$-operad $\PH^n .$
 The object $a_T$ of the canonical internal operad  in $\PH_T^n$ is the terminal object in this category.
 The nerve of $\PH^n$ is obtained by a bar-construction on the terminal pruned $n$-operad. \end{corol}
\noindent A more explicit description of $\PH^n$ will be given later. 

The functor of desymmetrisation for general  $(n-1)$-terminal $n$-operads 
  factorises through the category of pruned $(n-1)$-terminal $n$-operads
\begin{equation}\label{padj1} SO(V) \stackrel{Des_n}{\longrightarrow} \PO(V) \hookrightarrow O^{(n-1)}_n(V)\end{equation}
where $Des_n$ is defined by the formulas identical to the formulas for the desymmetrisation functor from  \cite{BEH}.  
 If $V$ is cocomplete the inclusion $\PO(V) \hookrightarrow O_n(V)$ has a left adjoint  $L.$ 
 It follows that the symmetrisation  functor can also  be factorised  
\begin{equation}\label{padj2} SO(V) \stackrel{Sym_n}{\longleftarrow} \PO(V) \stackrel{L}{\longleftarrow} O_n(V)\end{equation}
Notice that we  use the same notations for pruned versions of symmetrisation and desymetrisation functors 
as we used for the unpruned case in \cite{BEH}. We believe that this does not lead 
 to  confusion since $\PO(V)$ is a full subcategory of $O^{(n-1)}_n(V) .$ 
 Moreover, we will use the same notation in the reduced case  in Section \ref{reducedoperads}.

Again, as in the upruned case \cite{BEH}, we get the 
  following commutative square of adjunctions:

  \begin{center} {\unitlength=1mm
\begin{picture}(60,26)
\put(10,20){\makebox(0,0){\mbox{$SO(Set) $}}}
\put(9,6){\vector(0,1){10}}
\put(11,16){\vector(0,-1){10}}

\put(45,20){\makebox(0,0){\mbox{$\PO(Set)$}}}
\put(42,6){\vector(0,1){10}}
\put(44,16){\vector(0,-1){10}}

\put(18,21){\vector(1,0){14}}
\put(32,19){\vector(-1,0){14}}
\put(22,22){\shortstack{\mbox{$\scriptstyle Des_n $}}}
\put(22,16.7){\shortstack{\mbox{$\scriptstyle Sym_n $}}}

\put(68,10){\shortstack{\mbox{$ $}}}
\put(35,10){\shortstack{\mbox{$ $}}}
\put(10,3){\makebox(0,0){\mbox{$Coll_1(Set)$}}}
\put(45,3){\makebox(0,0){\mbox{$PColl_n(Set)$}}}
\put(19.5,4.5){\vector(1,0){14}}
\put(33.5,2.2){\vector(-1,0){14}}
\put(24,5.4){\shortstack{\mbox{$\scriptstyle W_n $}}}
\put(24,0.1){\shortstack{\mbox{$\scriptstyle C_n $}}}

\put(4,10){\shortstack{\mbox{$\scriptstyle {\cal F}_{\infty} $}}}
\put(12,10){\shortstack{\mbox{$\scriptstyle U_{\infty} $}}}

\put(35,10){\shortstack{\mbox{$\scriptstyle {\cal PF}_n $}}}
\put(45,10){\shortstack{$
{\scriptstyle PU_n} \addtocounter{equation}{1}
 \ \ \ \ \ \ \ \ \ \ \ \ \ \ \ \  \ \ \ \ \ \ \ \ \  \ \  (\theequation)
$}}

\end{picture}}
\end{center}\label{squareone} 

\noindent So, again, by the general theory of internal algebras,  we have a representable $2$-functor 
of  internal pruned $n$-operads on the $2$-category of  symmetric $Cat$-operads \cite{BEH}.  
We will denote by $\ph^n$ the symmetric categorical operad which  represents this $2$-functor. 

The value of the   composite $\F_{\infty}C_n$ on a pruned $n$-collection $X$ is easy to describe. In  arity $k$ 
it consists of all labelled planar trees with label $\{1,\ldots,k\}$ decorated by pruned $n$-trees and elements of $X_T .$
So Theorem 2.3 from \cite{BEH} provides us with  the following  description of   the operad $\ph^n .$ 

The objects of $\ph^n$ are labelled planar trees decorated by pruned $n$-trees. Analogously to the categorical operad  $\h^n$ in 
\cite{BEH}, the morphisms in $\ph^n$ are generated by simultaneous contractions of the input edges of a vertex provided there exists a corresponding morphism in $\Omega_n .$
We also can grow an internal edge by introducing a decoration by the linear tree $U_n $ (see \cite{BEH} 
for the analogous description of $\h^n$).  

The description we gave for  the $n$-operad $\PH^n$ in Corollary \ref{PH} is not very revealing.  We are going to make it  
more accessible by using some more structured planar trees we will call {\it composable}.

Recall from \cite{BEH} that the square  \ref{prunedoperads}.\ref{squareone} induces a map 
 $$\beta: Obj(\PH^n) = \PF_n(1)\rightarrow  Des_n(\F_{\infty}(C_n (1)  )).$$
This is not an injective transformation but it shows that to every object 
of  $\PH^n$  (an admissible expression) we can naturally associate a planar tree decorated by $n$-ordinals. 
We will call it  the underlying tree of  the admissible expression.
  
We define a notion of composable tree by induction.
 We call a corolla decorated by an  $n$-tree  $T$ (or $e$)  a decorated tree {\it  composable to $T$}.
 Suppose we  have already defined a notion of  composable decorated tree for  which their underlying 
 trees have length less than or equal to $l .$ Let $\tau$ be an equivalence class of admisssible 
expressions whose underlying  decorated tree  has  length $l+1 .$ 
Suppose, that its root vertex $v_0$ is decorated by a pruned tree $S .$ We will call the tree {\it composable}
 to an $n$-tree $T$ if it is equipped with  a morphism of $n$-trees 
$\sigma: T \rightarrow S$ with fibers $T_i , \ 1\le i \le k ,$ and  for each $i$ the $i$-th  
branch  at $v_0$ is composable to $T^{(p)}_i .$
 
 \
  
 \Example 

\vskip10pt \hskip40pt \includegraphics[width=190pt]{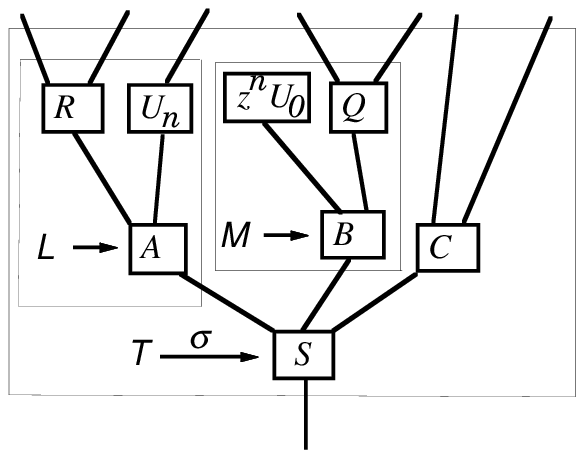}

\

The following lemma  is an $n$-operadic  analogue of the trivial fact that, given a  string of elements in a monoid, 
 there is a canonical way to calculate its value by performing multiplication always starting from  
the most right pair of elements. The value, however, does not depend on the method of multiplication. 
 We equally could choose as  canonical,   multiplication from the left end of the string. 
In the $n$-operad case, however, the situation is more subtle as  Example \ref{Tamarkintree}.2  shows.
 \begin{lemma}\label{composable} For every equivalence class of admissible expressions of arity $T ,$ there is a unique
 decorated tree composable  to $T$ representing this expression. \end{lemma} 
 
 \Proof  The proof is by a routine induction.  
 
 \
 
 \Q

\
 
 So one can think of the objects of $\PH^n_T $ as trees composable to $T .$ The morphisms are generated 
by composition of some nodes but in contrast with symmetric operads this operation can give a 
 tree which lies outside $\PH^n_T .$   

\
 
\Example \label{Tamarkintree} The following  tree  is composable but we cannot produce the same composition if we start to compose it from its root to its leaves. 
 This is actually the  combinatorial `raison d'\^{e}tre' of  Tamarkin's counterexample to cellularity of the Getzler-Jones operad.  

\
 
 \vskip15pt \hskip40pt \includegraphics[width=210pt]{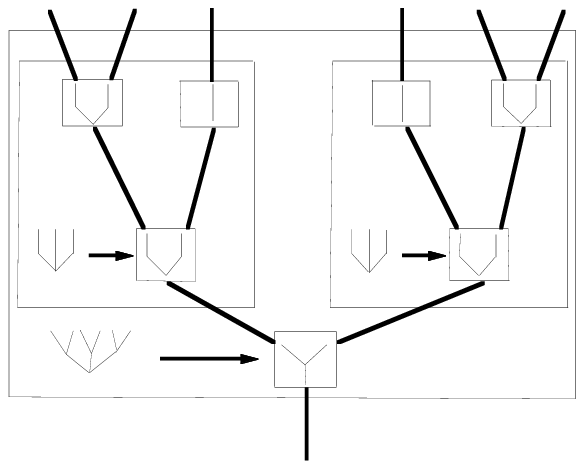}

 \

We also formulate some results about $\ph^n $ which can be proved following verbatim the proofs of the parallel results for 
  $\h^n$ in  \cite{BEH}:

 \begin{theorem} There is a natural isomorphism $Sym_n(\PH^n) \simeq \ph^n .$
 This isomorphism induces an isomorphism of  nerves
   $$N(\ph^n) \rightarrow Sym_n(N(\PH^n)) .$$
\end{theorem}
  \begin{theorem}\label{formulap} Let $A$ be a cocomplete symmetric $Cat$-operad and 
$a$ be an internal pruned 
$n$-operad in $A .$ Then 
$$(Sym_n(a))_k \ \simeq \ \mbox{\rm co}\hspace{-0.5mm}\!\lim\limits_{\ph^{n}_k} \tilde{a}_k$$
where $\tilde{a}_k:\ph^{n}_k \rightarrow A_k$ is the operadic functor generated by the  operad $a$.
\end{theorem}

Finally,  coming back to the adjunctions (\ref{padj1}),(\ref{padj2})  we find  one more interesting categorical 
pruned $n$-operad in this picture, namely $L(\H^n) .$ It comes together with  a canonical morphism $L(\H^n) \rightarrow \PH^n .$
 Obviously, $L(\H^n)$ is a categorical pruned $n$-operad freely generated by an internal (not pruned !)
  $n$-operad. 
 From this characterisation 
 we deduce a description of $L(\H^n).$ The objects are elements of the  free pruned $n$-operad on the 
following pruned $n$-collection: for a pruned tree $T$ it consists of all full injections $T\rightarrow  S .$ 
So, for any nonpruned $n$-tree $S ,$ we have an object, which we will denote by $S$ as well,  in the category $L(\H^n)_{S^{(p)}} $
 and these objects form an internal $n$-operad. The generators for morphisms  correspond to the morphisms of nonpruned trees and have the form
 $\mu(S; T_1,\ldots,T_p) \rightarrow T$ for a morphism of $n$-trees $\sigma: T\rightarrow S .$
 Notice, however, that this morphism generates a morphism between admissible expressions  
which corresponds to $\sigma^{(p)}: T^{(p)}\rightarrow S^{(p)} .$
 
\begin{conj} The operad  $N(L(\H^n))$  is a contractible simplicial operad. \end{conj}
 
 If this conjecture is true then we would be able to prove that the operad $N(\h^n)$ of \cite{BEH} 
is equivalent to the little $n$-disks operad. All our efforts to prove this conjecture have failed so far 
because some morphisms in $L(\H^n)$ are going  in opposite directions, which creates a 
 lot of combinatorial difficulties in analysing the homotopy type of $L(\H^n) .$         

\section{Reduced $(n-1)$-terminal  $n$-operads }\label{reducedoperads}

 We will call a symmetric operad  $A$  {\it reduced} if $A_0=A_1 = I $ the unit for tensor product, and the operadic unit is given by the identity.  
 Equivalently a reduced symmetric operad can be described  as a contravariant functor $A$ from the subcategory 
of nonempty ordinals and injective morphisms of the category $\Omega^s$ \cite{BEH}  such that $A([1]) = I ,$
 ({\it a reduced symmetric collection}) plus operadic composition for any surjection of finite sets.  
Operadic composition has to be natural  with respect to the injections.  The maps of reduced symmetric operads 
(reduced symmetric collections) are the maps of operads (natural transformations) which induce  identity morphisms in arity 
 $0$ and $1 .$    
 We use the notation $RSO(V)$ for the category of reduced symmetric operads and  $RColl_{\infty}(V)$
 for the category of reduced symmetric collections.
 
Observe that our category of reduced operads is just a subcategory of the category of reduced operads of Berger and Moerdijk  
\cite{Berger,BergerM} since they do not require $A_1$ to be a unit of $V .$ 

We also have a category of  reduced nonsymmetric collections which we will denote by $RColl_1(V) .$ 
These are contravariant functors on injective maps of $\Delta^{+}$ with the conditions $A([1]) = I .$
Recall from  \cite{ML}, that $\Delta^{+}$ is the full subcategory of $\Delta$ with nonempty ordinals as objects.

The $n$-operadic counterpart  of $RSO(V)$ will be the  category $\RO(V)$ of reduced   $(n-1)$-terminal $n$-operads.

\begin{defin} A pruned  $(n-1)$-terminal  $n$-operad $A$ is called  reduced   if 
$$A_{z^n U_0} = A_{U_n} = I $$ 
and its unit is given by the identity.
A morphism between two reduced $n$-operads is an $n$-operadic  morphism which induces  identity morphisms in arity  
$z^n U_0$ and $U_n .$ \end{defin}

\begin{defin} Let $\Lambda_n^{inj}$ be the category of  pruned nondegenerate $n$-trees and their  injective morphisms. 
A reduced $n$-collection is a contravariant functor $A$ from $\Lambda_n^{inj}$ to $V $  such that 
$A_{U_n} = I .$
\end{defin}  

The definition of the category of {\it reduced $n$-collections} $RCall_n(V)$ is obvious.

As for the reduced symmetric operads we have an equivalent characterisation of the reduced $(n-1)$-terminal operads.

\begin{pro} A reduced $(n-1)$-terminal $n$-operad   is given by a reduced $n$-collection $A$
together with a   multiplication 
$$\mu_{\sigma}: A_S \otimes A_{T_1^{(p)}}\otimes \ldots \otimes A_{T_k^{(p)}} \rightarrow A_T $$
for every surjection of $n$-trees $\sigma:T\rightarrow S ,$ satisfying associativity and unitary axioms 
and naturality with respect to injections of pruned trees. 
 \end{pro} 
  \begin{pro}\label{ron} 
  If $V$ is a cocomplete symmetric monoidal category then the forgetful functor $$RU_n:\RO(V) \rightarrow RColl_n(V)$$ has a  left adjoint $\RF_n$  and is monadic. 

          The free reduced $n$-operad monad on the category of reduced  $n$-collections in $Set$  is finitary and cartesian. 
  
 The analogous statements hold for symmetric operads and reduced nonsymmetric collections.  \end{pro}
 \Proof  Everything goes through  in full analogy with the unpruned case. 
The only difference is that we use only surjections of  pruned trees and we require some more identifications: 
 $$e \sim U_n .$$
 This identification leads to the effect that we do not have an underlying $n$-planar 
 tree with a vertex decorated by the linear tree $U_n.$

 \
 
 \Q
 
 \

 We also have a reduced desymmetrisation  functor  from the category  $RSO(V)$ to the category   
 $\RO(V)$ which we again will  denote by  $Des_n $ and 
 from general considerations of \cite{BEH}, we have a commutative square of adjunctions
  \begin{center} {\unitlength=1mm
\begin{picture}(60,26)
\put(10,20){\makebox(0,0){\mbox{$RSO(Set) $}}}
\put(9,6){\vector(0,1){10}}
\put(11,16){\vector(0,-1){10}}

\put(45,20){\makebox(0,0){\mbox{$\RO(Set)$}}}
\put(42,6){\vector(0,1){10}}
\put(44,16){\vector(0,-1){10}}

\put(19,21){\vector(1,0){14}}
\put(33,19){\vector(-1,0){14}}
\put(22,22){\shortstack{\mbox{$\scriptstyle Des_n $}}}
\put(22,16.7){\shortstack{\mbox{$\scriptstyle Sym_n $}}}

\put(35,10){\shortstack{\mbox{$ $}}}
\put(10,3){\makebox(0,0){\mbox{$RColl_1(Set)$}}}
\put(45,3){\makebox(0,0){\mbox{$RColl_n(Set)$}}}
\put(57,2){\shortstack{\mbox{$. $}}}
\put(20.5,4.5){\vector(1,0){14}}
\put(34.5,2.2){\vector(-1,0){14}}
\put(24,5.4){\shortstack{\mbox{$\scriptstyle W_n $}}}
\put(24,0.1){\shortstack{\mbox{$\scriptstyle C_n $}}}

\put(1,10){\shortstack{\mbox{$\scriptstyle {\cal RF}_{\infty} $}}}
\put(12,10){\shortstack{\mbox{$\scriptstyle RU_{\infty} $}}}

\put(35,10){\shortstack{\mbox{$\scriptstyle {\cal RF}_n $}}}
\put(45,10){\shortstack{\mbox{${\scriptstyle RU_n} 
\refstepcounter{equation}{1}
 \ \ \ \ \ \ \ \ \ \ \ \ \ \ \ \  \ \ \ \ \ \ \ \ \  \ \ (\theequation\label{radjsq}) $}}}

\end{picture}}
\end{center}

 \noindent Hence, we can develop the theory of internal reduced $n$-operads.  We denote by $\RH^n$ the reduced categorical $n$-operad freely generated by an internal reduced $n$-operad and by $\rh^n$ the symmetric categorical operad freely generated by an internal reduced $n$-operad.  

\begin{defin} An $n$-tree is called reduced if it  is pruned, nondegenerate and is not equal to $U_n .$ A planar tree decorated by $n$-ordinals is called reduced if all the decorations are reduced. 
 \end{defin}
 
 As in the  pruned case, we have 
 a natural transformation
 $$\gamma: \RF_n \rightarrow  Des_n(\RF_{\infty}(C_n  ))$$
and   have a description of the free reduced $n$-operad functor in terms of composable  reduced planar trees.  There are, however, some changes in the definition of composable trees due to the fact that we 
have an identity $e = U_n .$ 

Again we define the composable trees by induction on the length of the underlying tree, but this time for $l=1$ we 
 call a tree composable to $T$ when it is  a corolla decorated by an  $n$-tree  $S$   together with a quasibijection
 $T\rightarrow S .$  Suppose we  have already defined the notion of  composable decorated tree for  which the underlying  trees
 have length less than or equal to $l .$ Let $\tau$ be an equivalence class of admisssible expressions whose 
underlying  decorated tree  has  length $l+1 .$ Suppose  that its root vertex $v_0$ is 
decorated by a pruned tree $S .$ We will call the tree {\it composable} to an $n$-tree $T$ if it is equipped with 
 a surjection of $n$-trees $\sigma: T \rightarrow S$ with fibers $T_i , \ 1\le i \le k ,$ and  for each $i$ the $i$-th  branch  at $v_0$ is composable to $T^{(p)}_i .$ Notice that a branch can be 
 empty if the corresponding fiber of $\sigma$ is the linear tree $U_n .$

\begin{lemma}\label{composabler} For every equivalence class of admissible expressions of arity $T$ there is 
a unique reduced  tree composable  to $T$ representing this expression. \end{lemma} 

A nice result is the following characterisation of the free reduced $n$-operad functor.

 \begin{theorem}\label{submonadr}  The natural transformation $\gamma$ is injective.  The image of the inclusion 
 $$(\RF_n)_T(1)\subset Des_n(\RF_{\infty}(C_n (1) )_T$$
 consists of  labelled reduced  trees  decorated by $n$-ordinals which are 
 dominated by the $n$-ordinal $T .$ 
 \end{theorem}

\Proof
It is obvious that  a labelled decorated tree from the image of $\gamma$ is dominated by $T .$   We are going to describe a 
procedure which reconstructs  a unique composable reduced decorated tree from a labelled reduced 
planar $n$-tree dominated by $T ,$ so our theorem will follow from Lemma \ref{composabler}. 
  
 We construct a composable tree by induction. Let a labelled decorated tree have length equal to $1 .$ 
 So, this is a labelled corolla decorated by an $n$-ordinal $S .$  Suppose we know that this tree is dominated
 by an $n$-ordinal $T .$ The last property means that  the identity map of labels can be extended to a quasibijection 
 $\sigma:T\rightarrow S $  by Corollary \ref{quasibijection}.  So we have  a composable tree  $\mu(S;e, \ldots, e)$ of length $1$.
  
 Suppose we have already constructed a composable decorated tree of length $l$  for every labelled reduced $n$-planar 
  tree of length $l$ dominated by $T .$   Let $\tau$ be a labelled $n$-planar tree of length $l+1$ such that 
at the root $v_0 $  we have  an $n$-ordinal $S .$   Let us define a surjection  $f:|T|\rightarrow |S|$ as follows:
 $$f(i) = \#_{v_0}(i) .$$
The condition  $\tau \dom T$  implies  that  if $i<_p j$ in $T$ then  either $f_n(i)= f_n(j)$ or 
if they are not equal  then $$f(i) =  \#_{v_0}(i) <_r  \#_{v_0}(j) = f(j)$$ for $r\ge p$ or 
  $$f(j) =  \#_{v_0}(j) <_r  \#_{v_0}(i) = f(i)$$ for $r> p .$ So $f$ is a surjection  of ordinals. 
  
  Let $T_1, \ldots, T_k$ be  the list of fibers of $f .$ It is not hard to see  that the $i$-th branch of $\tau$ at 
the vertex $v_0$ is dominated by the $n$-ordinal $T_i .$  So we can proceed by induction and   finish the proof. 
   
   \
   
\Q

\

We now give a description of $\rh^n .$
The objects of $\rh^n$ are labelled planar trees decorated by reduced $n$-trees. 
 The morphisms are generated by simultaneous contractions of the input edges of
 a vertex provided there exists a corresponding surjection in $\Omega_n .$  
From this description we immediately get the following
\begin{lemma} The categories $\rh^n_k$ and $\RH^n_T$ are finite. \end{lemma}

There are no morphisms for growing  internal edge like  in the $\h^n$ and $\ph^n$ cases. 
For this reason the operad $\rh^n$ is even a finite poset operad but we do not need this property here. 
This is quite an important property, however, and we are going to  consider it in a separate paper. 

Similarly to the pruned case we have the following.

\begin{pro} The category $\RH^n_T$ is isomorphic to the comma-category of $\rh^n$ over its internal operad object 
$a_T .$  The $n$-operad structure  is given by the following 'convolution' product :

\noindent given a surjection of pruned trees $\sigma: T\rightarrow S$ and objects in comma-categories
$$x \rightarrow a_S \ , \  x_{1} \rightarrow a_{T_1^{(p)}} \ , \  \ldots \ , \  x_{k} \rightarrow a_{T^{(p)}_k} ,$$
we define an object of the  comma-category over $a_T$ by the composite:
$$\pi(\sigma)^{-1}m(x ; x_{1},\ldots,x_{k}) \rightarrow \pi(\sigma)^{-1}m(a_S;a_{T^{(p)}_1},\ldots,a_{T^{(p)}_k}) \rightarrow a_T$$
where $m$ is multiplication in $\rh^n$ and the last morphism is the structure morphism of the  internal operad $a$
in $\rh^n .$ 
       \end{pro}

 \begin{theorem}\label{Symnerve} There is a natural isomorphism $Sym_n(\RH^n) \simeq \rh^n .$
 This isomorphism induces an isomorphism of  nerves
   $$N(\rh^n) \rightarrow Sym_n(N(\RH^n)) .$$
\end{theorem}

\

We also want to  introduce an unbased version of reduced operads (we follow the terminology of \cite{Salvatore}). 
 These are reduced symmetric operads without nullary operations and  with $A_{1} = I .$ 
Notice, that we do not require $A_0$ to be $\emptyset$ but simply forget about the  $0$-arity of our operads. 
This is, of course, a linguistic difference but it  helps  to express the results nicely.     
\begin{defin} An unbased  reduced $n$-operad in a symmetric monoidal category $V$ consists of a
collection of objects $A_T\in V$, one for every reduced $n$-tree $T ,$  and $A_{U_n}=I ,$ together with a multiplication  
$$\mu_{\sigma}: A_S \otimes A_{T_1^{(p)}}\otimes \ldots \otimes A_{T_k^{(p)}} \rightarrow A_T $$
for every surjection of $n$-trees $\sigma:T\rightarrow S ,$ satisfying the usual associativity axiom and the unitarity 
 axiom with respect to the identity morphism $I \rightarrow A_{U_n} .$
\end{defin} 

We will denote the category of unbased reduced symmetric operads by $URSO(V)$ and the category of unbased reduced $n$-operads
 by $\URO(V).$ All the previous results about reduced operads can be carried over to the unbased case. 
In particular, we have a categorical unbased reduced operad 
 $\URH$ representing internal unbased $n$-operads and  a categorical unbased reduced symmetric operad $\urh.$ It is not hard to see, however, that there are canonical operadic maps
 $$\URH(+) \rightarrow \RH^n$$ 
 and
$$\urh(+) \rightarrow \rh^n ,$$
where $\URH(+)$ is obtained from $\URH$ by adding $\emptyset$ in the arity $z^nU_0$
($\urh(+)$ is obtained from $\urh$ by adding $\emptyset$ in the arity $0 $) 
and these maps are isomorphisms in  nondegenerate arities. 

\

Finally, we want to produce a version of the symmetrisation formula. This formula for reduced operads  admits a nice enhancement.

The opposite of the  Milgram poset $(\J^n_k)^{op}$ can be considered as  a subcategory of $\rh^n_k$ which consists
of labelled planar trees  with only one decoration. 
\begin{lemma}\label{cofinal}  For every $\tau\in \rh^n_k $
the comma-category $\tau/ (\J^n_k)^{op}$ is nonempty and connected. So
 $(\J^n_k)^{op}$ is a final subcategory
of $\rh^n_k$ \cite{ML}.
\end{lemma}

\Proof For $n=1$ the lemma is obviously true so  we can assume that $n\ge
2 .$

 It is not hard to see that  every object from $\rh^n_k$ is dominated by
one of the objects of $(\J^n_k)^{op} .$
An object $\tau\in \rh^n_k$ is a labelled planar tree decorated by
$n$-ordinals and, hence, determines
a canonical linear order on the set of its labels $\{1,\ldots,k\}.$
Without loss of generality we can assume that this ordered set is the
ordinal $[k] .$ So it is dominated by the $n$-ordinal
$M_0^k$ (see Section \ref{Comporders} for notation).

Let  $\tau \in \rh^n_k$ be dominated by $T'$ and $T'' .$

 Again without losing generality we  can assume that
$T'' = M^k_0$ and $T'$ can be obtained from $T''$ by permuting the labels.

Let us construct a totally 
$n$-ordered set $T$ which dominates $\tau $ and is dominated by both $ T'
$ and $T'' .$

To do this we apply the following reconstruction algorithm to $T' .$ Let
$i$ be  the first label
in $T'$ with respect to the linear order $<_0 $ and $j$ be the second. If
$i<j$ in $[k]$ then we define
$i<_0 j$ in $T^{(1)} .$
 If, however, $j<i$ in $[k]$ then we put $i<_1 j$ in $T^{(1)} .$ We also
put  all
 the other labels in $T^{(1)}$ to be  in the same order as in $T'.$ So we have 
constructed
an object of $(\J^n_k)^{op}$ which is dominated by $T'$   and  also dominates
$\tau .$

Now we continue this process and  do the same thing
for the second and third consecutive labels in $T'$ then for the third and
fourth and so on. In this way we construct
a sequence of $n$-ordered sets
$$\tau \dom  T =T^{(k-1)}\dom T^{(k-2)}\dom \ \ldots \  \dom T^{(1)}   \dom \ T^{(0)} = T'.$$
Finally, we observe, that $i<j$ implies $i<_0 j$ or $i<_1 j$ or $j<_1 i$
in $T$ by construction. This means
that $T\dom T''$ and we have finished the proof.

\

\Q

\

 \begin{theorem}\label{formular} Let $A$ be a cocomplete reduced symmetric  $Cat$-operad  and
$a$ be an internal reduced 
$n$-operad in $A .$ Then 
$$Sym_n(a)_k \ \simeq \ \mbox{\rm co}\hspace{-0.5mm}\!\lim\limits_{\rh^{n}_k} \tilde{a}_k
   \ \simeq \ \mbox{\rm co}\hspace{-2.4mm}\!\lim\limits_{\mbox{\rm (}\J^{n}_k\mbox{\rm )}^{\raisebox{0.1mm}{\scriptsize op}}} \tilde{a}_k                                                                     $$
where $\tilde{a}_k:\rh^{n}_k \rightarrow A_k$ is the operadic functor generated by  $a$.

The analogous formula holds in the  unbased case. \end{theorem}

To be able to apply this theorem to the reduced operads in a  symmetric  monoidal category $V,$ we have to exhibit  $V$ 
as a categorical reduced symmetric operad similar  to what was done in \cite{BEH} for the unreduced case. Our construction $V^{\bullet}$ 
 from \cite{BEH} produces only an unreduced operad. In the reduced  case we define the reduced categorical symmetric operad $V^{\bullet\bullet}$ as follows  
 \[
V^{\bullet\bullet}_k =  \left\{ 
\begin{array}{ll}
 V & \mbox{if} \  k\ge 2   \\
1 & \mbox{if} \  k=0,1 \end{array}
\right.
\]
with multiplication defined by iterated tensor product and trivial symmetric group action. It is obvious that there is an isomorphism of the categories of internal reduced operads in $V^{\bullet\bullet}$ and reduced 
operads in $V$ (see \cite{BEH} for this property in the unreduced case).  

\

A comment has to be made about endomorphism operads in the reduced situation. 
First of all we have to consider only the pointed case; i.e. we consider an object $X$ of our closed symmetric monoidal category 
 $V$ together with a fixed morphism $I\rightarrow X .$ For every $k\ge 0$ one can consider the following pullback 

\begin{center} {\unitlength=1mm
\begin{picture}(60,26)
\put(8,20){\makebox(0,0){\mbox{$\bar{V}(X^k,X) $}}}

\put(8,16){\vector(0,-1){10}}

\put(40,20){\makebox(0,0){\mbox{$V(X^k,X)$}}}

\put(40,16){\vector(0,-1){10}}

\put(17,20){\vector(1,0){14}}

\put(22,22){\shortstack{\mbox{$ $}}}

\put(68,10){\shortstack{\mbox{$ $}}}
\put(35,10){\shortstack{\mbox{$ $}}}
\put(8,3){\makebox(0,0){\mbox{$V(I^k,I)$}}}
\put(40,3){\makebox(0,0){\mbox{$V(I^k,X)$}}}
\put(17,3){\vector(1,0){14}}

\put(24,5.4){\shortstack{\mbox{$ $}}}

\put(13,10){\shortstack{\mbox{$ $}}}

\put(40,10){\shortstack{\mbox{$ $}}}

\end{picture}}
\end{center}
\noindent This symmetric collection has an obvious structure of an operad and we define 
{\it the  reduced endomorphism symmetric operad} $REnd(X)$ as follows:
 $$REnd(X)_k = \bar{V}(X^k,X) \ \mbox{for} \  k\ne 1  $$
 and  $$REnd(X)_1 = I .$$
 An action of a reduced operad $A$ on a pointed object $X$ is then 
 an operadic map $A\rightarrow REnd(X) .$ Usually an action of an operad   on a pointed object 
 is defined as an operad map from $A$ to its full endomorphism operad.  
 It is obvious, however, that such an action can be factorised through $REnd(X),$ so our definition agrees with the usual one.

 We give the analogous obvious definition for the reduced endomorphism $n$-operad $REnd_n(X).$ 
 Now, we will have a canonical isomorphism $$Des_n( REnd(X))\simeq REnd_n(X)$$ as in the unreduced case.  And, hence,  similarly to the unreduced case we have 
 \begin{theorem} For a reduced $(n-1)$-terminal  $n$-operad $A ,$ the categories of $Sym_n(A)$-algebras and  $A$-algebras are isomorphic. \end{theorem}

\section{Model structures on various categories of operads } 

Here we adapt the theory of \cite{BergerM} to the case of  $n$-operads.
Recall that one of the main technical points of this theory is the existence of a cofibrantly
 generated model structure on the category of reduced symmetric operads in a monoidal model category $V .$ 
This structure is transferred along the free symmetric operad functor from a  model structure  on symmetric collections in $V $ 
which in its turn is transferred along the free symmetric collection functor from the category of nonsymmetric collections. 
The term transferred means that the weak equivalences (fibrations) are precisely the morphisms which become weak equivalences 
(fibrations) after application of the right adjoint functor.

\begin{theorem}[Theorem 3.1 \cite{BergerM}]\label{BM1} There exists a transferred  model structure on the category of reduced symmetric operads in a monoidal model category $V $  if 
 \begin{itemize}
\item $V$ is cofibrantly generated and its unit $I$ is cofibrant;
\item   The comma category $V/I$ has a symmetric monoidal fibrant replacement functor;
\item $V$ admits a commutative Hopf interval.
\end{itemize}
\end{theorem}
\begin{theorem}[Theorem 3.2 \cite{BergerM}]\label{BM2}
If $V$ is a cartesian closed model category then there is a transferred model structure on the category of all symmetric operads in $V$ provided \begin{itemize}
\item $V$ is cofibrantly generated and the terminal object of $V$ is cofibrant;
\item $V$ has a symmetric monoidal fibrant replacement functor.
\end{itemize}
\end{theorem}

We already pointed out that our notion of reduced symmetric operad is stronger than that of \cite{BergerM}.
It is, however, not too hard to check that their proof works well in our situation. We only have to take care of
 their construction of a cooperad $TH$ from a commutative Hopf object $H$ \cite[Proposition 1.1]{BergerM}. 
 If we put $T'H_1 = I , T'H_n = H^{\otimes n} $ then $T'H$ is still a cooperad and it works for the reduced operads 
in our sense like $TH$ works in \cite{BergerM}.   So  Theorem \ref{BM1} holds in our case without any changes. 

Moreover, following  the  method of \cite{BergerM} in the case  of $n$-operads we can  prove:
\begin{theorem}\label{BMn1} There exists a transferred  model structure on the category of reduced $(n-1)$-terminal 
$n$-operads in a monoidal model category $V $  if 
 \begin{itemize}
\item $V$ is cofibrantly generated and its unit $I$ is cofibrant;
\item   The comma category $V/I$ has a symmetric monoidal fibrant replacement functor;
\item $V$ admits a commutative Hopf interval.
\end{itemize}
Moreover, the commutative square   \ref{radjsq} is a square of Quillen adjunctions.

The analogous theorem holds in the unbased case. 
\end{theorem}

\begin{theorem}\label{BMn2}
If $V$ is a cartesian closed model category then there is a transferred model structure on the categories of   
$(n-1)$-terminal
$n$-operads in $V$ and pruned $(n-1)$-terminal $n$-operads in $V$ provided \begin{itemize}
\item $V$ is cofibrantly generated and the terminal object of $V$ is cofibrant;
\item $V$ has a symmetric monoidal fibrant replacement functor.
\end{itemize}
The corresponding commutative squares of adjunctions are squares of Quillen adjunctions.
\end{theorem}

 Examples of  monoidal categories satisfying the conditions of these theorems are given in \cite{BergerM}. 
 The most important for us are the categories of simplicial sets,  the category $Top$ of  compactly generated topological spaces,
and the category $Ch(R)  $ of chain complexes over a commutative ring $R$ with unit. 
 We  note that the category of chain complexes satisfies the hypothesis of Theorem
\ref{BM1} and the categories of simplicial sets and topological spaces satisfy the asssumptions of  Theorem \ref{BM2}. 

 The weak equivalences (fibrations) between topological or simplicial $n$-operads are therefore 
operadic maps which are termwise weak equivalences (fibrations) in  simplicial sets and $Top.$  
The weak equivalences (fibrations) between chain $n$-operads are those operadic maps which are termwise quasi-isomorphisms (epimorphisms).

\  

Since many of our categorical operads are given by a bar construction we would like to investigate now what are the 
homotopy theoretic properties of the bar-construction for $n$-operads in general. We will do it for the general case of $(n-1)$-terminal
 $n$-operads. The pruned and reduced cases are similar.

In addition, suppose  $V$   is enriched in simplicial sets with  simplicial $hom$-functor  $V^S(-,-) $ satisfying 
\begin{equation}\label{simplenrich}V^S(X,Y) = V^S(I , V(X,Y)) \end{equation}
where $V(X,Y)$ is the internal $hom$-functor in $V .$    Then the categories  of  $n$-collections and $n$-operads  
become simplicially enriched. 
For two $(n-1)$-terminal $n$-collections $X,Y ,$ we define its simplicial set of morphisms
as 
$$\Coll(X,Y) = \prod_{T\in Tr_n} V^S(X_T,Y_T).$$

As in \cite{BEH}, let $(\F_{n}, \mu, \epsilon)$ be the free $(n-1)$-terminal $n$-operad monad. 
Then for two $n$-operads $A,B$ we define their simplicial set of morphisms as the equalizer
\begin{center}{\unitlength=0,8mm

\begin{picture}(90,10)(0,15)

\put(78,20){\makebox(0,0){\mbox{$ \Coll({\F}_{n}A,B)$}}}

\put(38,20){\makebox(0,0){\mbox{$ \Coll(A,B)$}}}

\put(53,20.8){\vector(1,0){6}}

\put(53,19.2){\vector(1,0){6}}

\put(2,20){\makebox(0,0){\mbox{$ \Oper(A,B)$}}}

\put(17,20){\vector(1,0){6}}

\end{picture}}

\end{center}

\noindent with the obvious horizontal morphisms generated by the operadic structures of $A$ and $B$. 

If $V$ has tensors and cotensors with respect to simplicial sets then the categories $O_n^{(n-1)}(V)$ and $Coll_n^{(n-1)}(V)$ also have them.
In particular, one can speak about the total object for a simplicial $n$-operad ($n$-collection) $A^{\star}$. By definition this is 
the coend
$$Tot(A^{\star}) = \int^{[n]\in\Delta^{+}} A^{n}\otimes \Delta^{n} ,$$ 
where   $\Delta^{n} =
\Delta^{+}([n+1],-)$ is the standard simplex of dimension
$n$.

Let $X$ be an $(n-1)$-terminal $n$-operad in $C$. Then the bar-construction $B(\F_{n},\F_{n},X) $ of $X$ is the total object of the simplicial operad

{\unitlength=1mm

\begin{picture}(60,12)



\put(19,5){\makebox(0,0){\mbox{${\cal F}_{n}(X)$}}}

\put(35,5.5){\vector(-1,0){7}}
\put(35,4.5){\vector(-1,0){7}}

\put(36,6.5){\shortstack{\mbox{\small $ $}}}

\put(45,5){\makebox(0,0){\mbox{${\cal F}_{n}^2(X)$}}}

\put(71.5,5){\makebox(0,0){\mbox{${\cal F}_{n}^3(X)$}}}

\put(61.5,5.8){\vector(-1,0){7}}
\put(61.5,5){\vector(-1,0){7}}
\put(61.5,4.2){\vector(-1,0){7}}

\put(83,5){\shortstack{\mbox{\small $\ldots $}}}

\put(54,7){\shortstack{\mbox{\small $ $}}}

\end{picture}}

\begin{theorem}  Let $V$ be a  model category which satisfies the conditions of  Theorem \ref{BMn2},
 has a simplicial enrichment $V^S(-,-)$ satisfying \ref{simplenrich}, and which is a simplicial 
model category with respect to these structures \cite{H}.  
Let $X$ be an  $(n-1)$-terminal $n$-operad  in $V$ with  cofibrant underlying $n$-collection. Then 
the canonical operad morphism
$$\rho: B(\F_{n},\F_{n},X)\longrightarrow X$$
is  a  cofibrant replacement for $X$ in the model category of $(n-1)$-terminal $n$-operads.

The analogous theorem holds in the pruned case and in the reduced and unbased reduced case if $V$ satisfies 
the assumptions of  Theorem \ref{BMn1}.    
\end{theorem}

\Proof It follows from the general properties of bar-construction \cite{May}  that the  morphism
$\rho$ is a deformation retraction in $Coll^{(n-1)}_n$. Hence, $\rho$ is a trivial fibration of operads.

It remains  to prove that
$B(\F_{n},\F_{n},X)$ is a cofibrant
$n$-operad.  Let $f:E\rightarrow B$ be a trivial fibration of $n$-operads.

We have to show that  any operadic map $B(\F_{n},\F_{n},X)\rightarrow B$ can be lifted to $E .$

{\unitlength=1mm

\begin{picture}(60,30)

\put(69,7){\makebox(0,0){\mbox{$B$}}}
\put(69,22){\vector(0,-1){12}}

\put(69,25){\makebox(0,0){\mbox{$E $}}}

\put(55,7){\vector(1,0){10}}

\put(71,15){\shortstack{\mbox{$f $}}}

\put(35,7){\makebox(0,0){\mbox{$B(\F_{n},\F_{n},X)$}}}

\put(48,10){\line(3,2){5}}
\put(54,14){\line(3,2){5}}
\put(60,18){\vector(3,2){5}}

\end{picture}}

\noindent By construction this amounts to the following lifting problem in the
category of cosimplicial spaces

{\unitlength=1mm

\begin{picture}(60,30)

\put(72,7){\makebox(0,0){\mbox{$\Oper(\F_{n}^{\star}X,B)$}}}
\put(72,22){\vector(0,-1){12}}

\put(72,25){\makebox(0,0){\mbox{$\Oper(\F_{n}^{\star}X,E)$}}}

\put(43,7){\vector(1,0){12}}

\put(74,15){\shortstack{\mbox{$f^{\star} $}}}

\put(40,7){\makebox(0,0){\mbox{$\Delta^{\star}$}}}

\put(43,10){\line(3,2){5}}
\put(49,14){\line(3,2){5}}
\put(55,18){\vector(3,2){5}}

\end{picture}}
  
\noindent where $\Delta^{\star}$ is the cosimplicial simplicial set consisting of standard simplices.
Since  $\Delta^{\star}$ is cofibrant in the Reedy model structure \cite{H} it remains to
show that
$f^{\star}$ is a trivial fibration. 

 We follow  a method developed in \cite{BA}. We have
to prove that in the diagram

{\unitlength=1mm

\begin{picture}(100,45)(-10,0)

\put(72,7){\makebox(0,0){\mbox{$M_i(\Oper(\F_{n}^{\star}X,B))$}}}
\put(72,22){\vector(0,-1){12}}

\put(72,25){\makebox(0,0){\mbox{$M_i(\Oper(\F_{n}^{\star}X,E))$}}}

\put(41,7){\vector(1,0){11}}

\put(74,15){\shortstack{\mbox{$M_i f^{\star} $}}}

\put(8,35){\vector(1,-2){12.5}}
\put(22,35){\vector(4,-1){30}}
\put(14,35){\vector(1,-1){7}}

\put(25,7){\makebox(0,0){\mbox{$\Oper(\F_{n}^{i+1}X,B)$}}}

\put(12,38){\makebox(0,0){\mbox{$\Oper(\F_{n}^{i+1}X,E)$}}}

\put(25,25){\makebox(0,0){\mbox{$W_i$}}}
\put(25,22){\vector(0,-1){12}}
\put(30,25){\vector(1,0){22}}

\put(30.5,20){\line(-1,0){4}}
\put(30.5,20){\line(0,1){3.5}}

\put(22,31){\makebox(0,0){\mbox{$\omega_i$}}}

\end{picture}}

\noindent the canonical map $\omega_i$ to the pullback is a trivial
fibration. In this diagram $M_i(-)$ is the $i$-th  matching object of
the corresponding cosimplicial object \cite{H}. 

According to Lemma 2.3 from \cite{BA} the diagram above is isomorphic
to the diagram

{\unitlength=1mm

\begin{picture}(100,45)(-10,0)

\put(71,7){\makebox(0,0){\mbox{$\Coll(L_{i-1} \F_{n}^{\star -1}X,B)$}}}
\put(71,22){\vector(0,-1){12}}

\put(71,25){\makebox(0,0){\mbox{$\Coll(L_{i-1} \F_{n}^{\star -1}X,E)$}}}

\put(40,7){\vector(1,0){10}}

\put(73,15){\shortstack{\mbox{$M_i f^{\star} $}}}

\put(8,35){\vector(1,-2){12.5}}
\put(22,35){\vector(4,-1){28}}
\put(14,35){\vector(1,-1){7}}

\put(25,7){\makebox(0,0){\mbox{$\Coll(\F_{n}^{i}X,B)$}}}

\put(12,38){\makebox(0,0){\mbox{$\Coll(\F_{n}^{i}X,E)$}}}

\put(25,25){\makebox(0,0){\mbox{$W_i$}}}
\put(25,22){\vector(0,-1){12}}
\put(28,25){\vector(1,0){22}}

\put(30.5,20){\line(-1,0){4}}
\put(30.5,20){\line(0,1){3.5}}

\put(22,31){\makebox(0,0){\mbox{$\omega_i$}}}

\put(40,33){\makebox(0,0){\mbox{$\phi_i$}}}
\put(44,9){\makebox(0,0){\mbox{$\psi_i$}}}

\end{picture}}

\noindent Here,  $L_i(\F_{n}^{\star -1}(X))$ is the latching object \cite{H} for the augmented cosimplicial
object $\F_{n}^{\star -1}(X)$ in $\Coll(V) ,$
 and $\phi_i, \psi_i$ are generated by the canonical morphism
$$l_{i-1}: L_{i-1} \F_{n}^{\star-1}X\rightarrow \F_{n}^i X .$$
If we show that this morphism is a cofibration, then $\omega_i$ will be a trivial fibration by
the axioms for a simplicial model category. 

We will actually  prove that  $l_{i-1}$ is an isomorphism onto a summand.

 It was proved in \cite{BEH}[Theorem 9.1] that the $k$-th iteration of the  
functor $\F_{n}$ is given by the following formula:
\begin{equation}\label{freeop}
{\cal F}_{n}^k(X)_T = \coprod_{W_1\stackrel{{\scriptstyle f}_{1}}{\leftarrow\!-}W_2\stackrel{\scriptstyle f_2}{\leftarrow\!-}\
\ldots\stackrel{\scriptstyle f_{\scriptscriptstyle k-1}}{\leftarrow\!-}\ W_k}
\tilde{X}(W_k),\end{equation}
where $f_1,\ldots,f_{k-1}$ are morphisms in $\H^n_T .$

 The coface operators in $\F_{n}^{\star -1}(X)$ are canonical inclusions on the summands corresponding to the
operators of insertion of the identities to the chain 
 $$W_1\stackrel{{\scriptstyle f}_{1}}{\leftarrow\!\!\!\!-}W_2\stackrel{\scriptstyle f_2}{\leftarrow\!\!\!\!-}\
\ldots\stackrel{\scriptstyle f_{\scriptscriptstyle k-1}}{\leftarrow\!\!\!\!-}\ W_k.$$ 
The rest of the proof follows in  complete analogy with Lemma 4.1 of \cite{BA}.
  
 The proof in the reduced and unbased cases are analogous. Notice that, in our version of the category of reduced operads,
any reduced collection is automatically well pointed in the sense of \cite{BergerM} so we do not include this condition 
in the formulation of the reduced version of our theorem. 
 
 \
 
\Q

\

\begin{theorem} The simplicial $n$-operads $N(\H^n), N(\PH^n), N(\URH)$ are cofibrant $n$-operads in the categories 
of  $(n-1)$-terminal simplicial $n$-operads, $(n-1)$-terminal simplicial pruned $n$-operads and 
unbased reduced $(n-1)$-terminal simplicial pruned $n$-operads respectively. 

The simplicial symmetric operads $N(\h^n), N(\ph^n)$ are cofibrant simplicial symmetric operads. 

The unbased reduced symmetric operad $N(\urh)$ is cofibrant in the category of unbased reduced symmetric operads. 

The same theorem is true for the geometric realisations of these operads and for 
the  reduced operad of chain complexes of $\URH$ and $\urh .$ 
\end{theorem}

\Proof The nerves of the above $n$-operads are bar-constructions on the terminal $n$-collection in the corresponding categories. 
In these categories the terminal $n$-collection is obviously a cofibrant collection.

The cofibrantness of $N(\h^n), N(\ph^n)$ and $N(\urh)$ follows from the corresponding Quillen adjunctions between 
$n$-operads and symmetric operads and the fact that in these cases the symmetrisation functor commutes with nerves. 

The topological and chain versions of the theorem follow from the general considerations of \cite{BergerM}.

\

\Q

\

\Remark  The nerves of $\RH^n$ and $\rh^n$ are not cofibrant.  The reason is that $1$ is not a cofibrant reduced $n$-collection.

\section{Fulton-Macpherson operad and  Getzler-Jones decomposition }\label{FMandGJ}
 
 The operadic structure on compactified moduli space of configurations of points in $\Re^n$ was first observed by Getzler and Jones in \cite{GJ}.  
 Here we use an explicit approach of \cite{KS,Sinha} to describe this compactification.

Let $\mod{k}$ be  the  quotient of the configuration space
 $$\Conf_k(\Re^n)= \{(x_1,\ldots, x_k)\in (\Re^n)^k \ | \  x_i \ne x_j  \  \mbox{if} \  i\ne j \ \}$$
 with respect to the obvious action of the $(n+1)$-dimensional Lie group $G_n$ of affine transformations of the form $u \mapsto
 \lambda u + v ,$ where $\lambda > 0$ is a real number and $v$ is a vector from $\Re^n .$ 
 The functions 
 $$u_{ij}(x_1,\ldots,x_k) =  \frac{x_j - x_i}{||x_j-x_i||} \ \ \  \ \ \ \    \ 1\le i , j\le k \ , i\ne j \ \  \  \ \ \ \ \ \ \ \ \ \ \  \ \ \ \ $$  
 $$ d_{i,j,l}(x_1,\ldots,x_k) =  \frac{||x_i - x_j||}{||x_i-x_l||}  \  \ \ \ \ \ \ 1\le i,j,l\le k  \ , \  i\ne j , j\ne l , i \ne l  $$ 
 embed  $\mod{k}$ into a compact space
 $$\Gamma_k[\Re^n] = \prod_{\stackrel{\scriptstyle 1\le i,j \le k}{i\ne j}} S^{n-1}\times \prod_{\stackrel{\scriptstyle 1\le i,j,l\le k}{ i\ne j , j\ne l , i \ne l} }[0,\infty] . $$ 
The $k$-th space of the Fulton-Macpherson operad $\fm^n$ is obtained as the  closure
 of the moduli space $\mod{k}$ inside  $\Gamma_k[\Re^n].$ Notice, that we use here a
 reduced version of $\fm^n$ meaning that we put $\fm^n_0 = 1 .$ If we forget about 
nullary operations we will get the unbased version of the Fulton-Macpherson operad which we will denote by  $\fm^n_{\circ} .$  Notice also  that in \cite{GJ,KS} the Fulton-Macpherson operad means  $\fm^n_{\circ} .$

In \cite{GJ, KS, Salvatore, Sinha} the following  properties of $\fm^n$ were established 

\begin{itemize}
\item $\fm^n$ is a reduced symmetric  operad weakly equivalent to the little $n$-cube operad; 
\item $\fm^n_{\circ}$ is an unbased  reduced cofibrant operad weakly equivalent to the unbased little $n$-cube operad;
\item set theoretically $\fm^n$ is the free reduced operad  on the reduced symmetric collection $\mod{\bullet} ;$
\item $\fm^n_k$ is a manifold with corners;
\item $\fm^1$ is isomorphic to the Stasheff operad of associahedra.  
\end{itemize}

 Later in this paper we will use notations $u_{ij}$ and $d_{ijl}$ for the coordinates of points in $\fm^n .$
The following lemma describes the behaviour of these   coordinates  under operadic multiplication.   
\begin{lemma}\label{multformula} Let $\sigma:[n]\rightarrow [k]$ be a surjection  in $\Omega^s .$ Then 
 \[
u_{ij}(\mu(x; x_{1}, \ldots, x_{k})) =  \left\{ 
\begin{array}{ll}
  u_{ij}(x_{l})& \mbox{if} \ \sigma(i) = \sigma(j) = l    \\
u_{\sigma(i')\sigma(j')}(x_{}) & \mbox{if} \ \sigma(i) \ne \sigma(j) 
\end{array}
\right.
\]
where $i',j'$ are images of $i,j$ in  the fiber of $\sigma$ over $l .$

Similarly
 \[
d_{ijl}(\mu(x ; x_{1}, \ldots, x_{k})) =  \left\{ 
\begin{array}{ll}
  d_{i'j'l'}(x_{s})& \mbox{if} \ \sigma(i) = \sigma(j) = \sigma(l) = s    \\
  0 & \mbox{if} \ \sigma(i) = \sigma(j) \ne  \sigma(l)     \\
d_{\sigma(i)\sigma(j)\sigma(l)}(x_{}) & \mbox{if} \ \sigma(i) \ne \sigma(j)  \ne \sigma(l) \ \sigma(i) \ne \sigma(l)
\end{array}
\right.
\]
All other values of $d_{ijl}(\mu(x ; x_{1}, \ldots, x_{k}))$ can be deduced from the above table and the
 following relations between $d_{ijk} $ \cite{Sinha}: 
$$d_{ijk}d_{ikj} = d_{ijk}d_{ikl}d_{ilj} = d_{ijk}d_{jki}d_{kij} = 1 .$$

\end{lemma}

\Proof  The proof can be obtained  using the explicit formulas for the operadic multiplication in $\fm^n$ from 
\cite{Markl} or techniques from \cite{Sinha} .  

\

\Q
 
\ 
 
Following Joyal \cite{J} we  give a  definition of a generalised  $n$-tree. 
\begin{defin} A generalised $n$-tree $X$ is a chain of partially ordered sets and functions
$$ R^n \stackrel{\rho_{n-1}}{\longrightarrow} R^{n-1} \stackrel{\rho_{n-2}}{\longrightarrow} \ldots 
\stackrel{\rho_{1}}{\longrightarrow} R^1 \stackrel{\rho_{0}}{\longrightarrow} R^0 = 1 $$
such that the  induced order on $\rho_i^{-1}(a)$ is linear for all $0\le i \le n-1$ and $a\in R^{i} .$

A topological $n$-tree is an $n$-tree for which all $R^i$ are endowed with a topology and all $\rho_i$ are continuous functions.
\end{defin}

The definition of a morphism of generalised $n$-trees (also from \cite{J}) coincides verbatim with the definition of morphism of finite $n$-trees in \cite{BEH}.
 
 Consider the following topological $n$-tree : $\Re^{\le n}$
 $$ \Re^n \rightarrow \Re^{n-1} \rightarrow \ldots \rightarrow \Re^1 \rightarrow \Re^0 = 0 $$
 where morphisms are projections on the first coordinates.  We introduce the $n$-tree structure on $\Re^{\le n}$  by ordering 
the fiber of each projection according to its natural order. Now, for every pruned $n$-tree $T ,$ one can  
consider the space of all injective $n$-tree morphisms from  $T$ to $\Re^{\le n} .$ This space is the classical 
{\it Fox-Neuwirth cell} $\fn_T$  corresponding to $T .$
 We will consider it as an open submanifold of configurations of points with the 
labelling prescribed by the order in $T $ \cite{Berger,GJ,VoronovG}.      
 
 \
  
\Example 

\vskip20pt \hskip45pt \includegraphics[width=190pt]{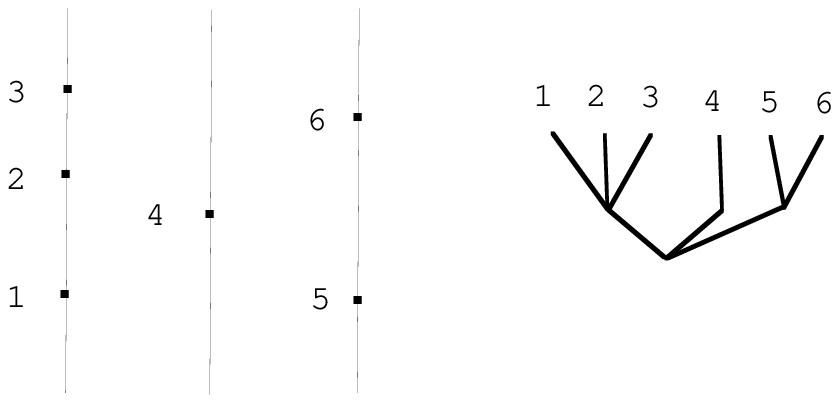}

\

The group $G_n$ obviously acts on $\fn_T$  and we call the corresponding quotient space {\it a Getzler-Jones cell }
 $\FN{T}.$  

Recall that the  functor $W_n: RColl_1 \rightarrow RColl_n , \ W_n(A)_T = A_{|T|}$ has a left adjoint 
$$C_n: RColl_n \rightarrow RColl_1 \ , \  C_n (B)_{[k]} =  \coprod_{|T|=[k]} B_T .$$

The configuration space $\Conf_k(\Re^n)$ admits the classical Fox-Neuwirth decomposition 
 $$\Conf_k(\Re^n)= \bigcup_{|T|=[k] , \pi \in \Sigma_k} \pi \fn_T.$$ 
This   means that 
there is a bijective continuous map of collections 
$$\varepsilon: S(C_n(\fn_{\bullet})) \rightarrow  \Conf_{\bullet}(\Re^n),$$  
where $S:RColl_1 \rightarrow RColl_{\infty}$ is the functor of symmetrisation of reduced nonsymmetric collections.  The inverse bijection 
 is not continuous.    
Certainly $\varepsilon$ is equivariant with respect to the action of $G_n$ and so we have a stratification of the moduli space
 of configurations
$$\mod{k}= \bigcup_{|T|=[k] , \pi \in \Sigma_k} \pi \FN{T}$$
and a corresponding continuous bijection which we will denote by the same letter
\begin{equation}\label{moddecomp}\varepsilon: S(C_n(\FN{\bullet})) \rightarrow  \mod{\bullet} . \end{equation} 
 The free reduced symmetric operad functor $\RF_{\infty}$ is  factorized as $\RF_s \cdot S ,$ where $\RF_s$ is the 
free reduced symmetric operad functor on reduced symmetric collections. So we have a bijective continuous map of operads
 \begin{equation}\label{GJdecomp} \RF_{\infty}(C_n(\FN{\bullet})) \simeq \RF_s(S(C_n(\FN{\bullet}))) \rightarrow  \RF_s(\mod{\bullet}) \rightarrow \fm^n
 \end{equation}

We would like to get a description  of $\FN{T}$  in terms of the functions $u_{ij} .$
Let $\stackrel{\scriptscriptstyle o \  \ \ \ \ \ \  }{S^{n-p-1}_{+}}$ denote the  open $(n-p-1)$-hemisphere in $\Re^n$ , $0\le p \le n-1$:
\[  \stackrel{\scriptscriptstyle o \  \ \ \ \ \ \  }{S^{n-p-1}_{+}} = \left\{(x_1, \ldots, x_n)\in \Re^n \left \vert \ \begin{array}[c]{ll}
  x_1^2+\ldots+x_n^2 = 1&   \\
   x_{p+1}> 0 \ \mbox{and} \ x_i = 0  \ \mbox{if} \  1\le i \le p & 
  \end{array} \right. \right\} \]
  and let 
 \[  \stackrel{\scriptscriptstyle o \  \ \ \ \ \ \  }{S^{n-p-1}_{-}} = \left\{(x_1, \ldots, x_n)\in \Re^n \left \vert \ \begin{array}[c]{ll}
  x_1^2+\ldots+x_n^2 = 1&   \\
   x_{p+1}< 0 \ \mbox{and} \ x_i = 0  \ \mbox{if} \  1\le i \le p & 
  \end{array} \right. \right\} \raisebox{-3mm}{ \ .} \]

The  closure of $\stackrel{\scriptscriptstyle o \  \ \ \ \ \ \  }{S^{n-p-1}_{+}}$ will be denoted $S^{n-p-1}_{+}$
 and the closure of $\stackrel{\scriptscriptstyle o \  \ \ \ \ \ \  }{S^{n-p-1}_{-}}$ will be denoted $S^{n-p-1}_{-} .$
Observe, that 
$$\stackrel{\scriptscriptstyle o \  \ \ \ \ \ \  }{S^{n-r-1}_{+}} \  \subset  \ {S^{n-p-1}_{+}} \  \mbox{for} \ r\ge p \ \ \  \mbox{and} \ \ \
\stackrel{\scriptscriptstyle o \  \ \ \ \ \ \  }{S^{n-r-1}_{-}} \ \subset  \ {S^{n-p-1}_{+}} \   \mbox{for} \  r>p  .$$

\begin{lemma}\label{GJcellT}  For a pruned $n$-tree $T ,\ |T|=k  , \ k\ge 0 ,$ the Getzler-Jones cell $\FN{T}$ is equal to the  set
\[   \left\{ x\in \mod{k}  \left | 
 \ \begin{array}[c]{ll}
  \ u_{ij}(x)  \in \stackrel{\scriptscriptstyle o \  \ \ \ \ \ \  }{S^{n-p-1}_{+}}&  \ \mbox{\rm if} \ i<_p j  \ \mbox{\rm in} \ T  \\
  &  \\
  \ u_{ij}(x)  \in \stackrel{\scriptscriptstyle o \  \ \ \ \ \ \  }{S^{n-p-1}_{-}}&  \ \mbox{\rm if} \ j<_p i  \ \mbox{\rm in} \ T  
  \end{array}\right.  \right\} \raisebox{-7mm}{ \ .}\]
 It is a contractible open manifold of dimension $E(T)-n-1 $ where $E(T)$ is the number of edges in the tree $T .$

\end{lemma}

\Proof Obvious from the definition of $\FN{T}$ . 

\

\Q

\

\begin{defin}The dimension of an  $n$-pruned tree $T\ne U_n$ is the integer 
$$dim(T) = E(T)-n-1. $$
We also put $dim(U_n)= 0. $ \end{defin} 
 
In virtue of (\ref{GJdecomp})  the decomposition (\ref{moddecomp})  can be extended to  a decomposition  of $\fm^n$ \cite{GJ} . Following \cite{VoronovG} we will call it 
{\it the Getzler-Jones decomposition}.  
The cells of this decomposition are indexed  by labelled reduced planar trees i.e. by the objects of $ \rh^n$  with vertices decorated by  points of   Getzler-Jones cells.  We will call such a space indexed by an object $\tau\in \rh^n$  {\it a generalised Getzler-Jones cell}   and will denote it $\FN{\tau} .$  Since the generalised Getzler-Jones cells do not intersect each other we have a correctly defined map
$$\tau:\fm^n \rightarrow \ \ \mbox{Obj}(\rh^n) .$$

\begin{lemma}\label{GJtau} Let $x\in \FN{\tau}$ then $$u_{ij}(x)\in \stackrel{\scriptscriptstyle o \  \ \ \ \ \ \  }{S^{n-p-1}_{+}}$$ if
$i<_p j$ in $\tau(x)$ or  $$u_{ij}(x)\in \stackrel{\scriptscriptstyle o \  \ \ \ \ \ \  }{S^{n-p-1}_{-}}$$ if
$j<_p i$ in $\tau(x) .$ \end{lemma}

\Proof  The proof is easily obtained by induction on the length of $\tau$ and Lemmas \ref{multformula}, \ref{GJcellT}. 

\

\Q

\

The Getzler-Jones decomposition is  not a cellular decomposition or stratification of $\fm^n $ since the boundary of the closure of a Getzler-Jones cell may not be equal to the union of low dimensional Getzler-Jones cells as  was first observed by Tamarkin (see  \cite{VoronovG} for a description of Tamarkin's counterexample).
 
 However, the following is true\footnote{I am grateful to Ezra Getzler and Sasha Voronov who explained this fact to me.} :
 
 \begin{pro}\label{FNclosure}  The closure of the  Getzler-Jones cell $K_T = cl(\FN{T})$ 
is a manifold with corners homeomorphic  to a ball of dimension $dim(T) .$  \end{pro}

\Proof  We have to use the original description of $\fm_k^n$ in terms of the iterated blow-up of $(\Re^n)^k$  along its fat diagonal
 \cite{GJ,Salvatore,VoronovG}.  Consider the closure of the
Getzler-Jones cell $\FN{T}$ in $(\Re^n)^k .$
 The intersection of this subspace with each of the diagonals is a manifold, so  restricting the Fulton-Macpherson 
blow up procedure we get a manifold with corners  homeomorphic to a ball since the blow up does not change the topological type of the manifold.  This manifold is homeomorphic to the closure of $\FN{T} $ in $\fm^n .$

\

\Q

 \section{Getzler-Jones operad}\label{GJoperad}
  
Let us  take  
$\RF_n(\FN{\bullet })$ to be  the free reduced $n$-operad generated by the reduced $n$-collection $\FN{\bullet}$. Then  we have   a canonical inclusion of $n$-operads
$$\gamma: \RF_n(\FN{\bullet})\longrightarrow RDes_n(\RF_{\infty}(C_n(\FN{\bullet}) )).$$ 
Consider the composite $\Phi$
$$\RF_n(\FN{\bullet})\stackrel{\gamma}{\longrightarrow} RDes_n(\RF_{\infty}(C_n(\FN{\bullet}) ) {\longrightarrow} $$ $$ \longrightarrow RDes_n(\RF_s(\mod{\bullet}))\longrightarrow  RDes_n(\fm^n) $$
which is an injective continuous map of $n$-operads.

\begin{defin}  The Getzler-Jones  $n$-operad  $\GJ^n$ is the image of  $\Phi .$ This is a reduced $(n-1)$-terminal $n$-operad.
\end{defin}
If we forget about operations with degenerate arity,  we obtain  {\it the unbased Getzler-Jones  $n$-operad } $\GJ^n_{\circ} .$

\begin{pro} Let $T\ne U_n$ be a nondegenerate pruned $n$-tree.  The following topological spaces  are equal:
\begin{itemize}
\item $\GJ^n_T \ ;$
\item $\bigcup\limits_{\tau \hspace{0.4mm}  \dom \hspace{0.4mm}T} \FN{ \tau}  \ ;$
\item Kontsevich-Soibelman space \cite{KS}
\[ X_T =  \left\{ x\in \fm^n_{|T|}  \left \vert \ \begin{array}[c]{ll}
  \ u_{ij}(x)  \in {S^{n-p-1}_{+}}&  \ \mbox{\rm if} \ i<_p j  \ \mbox{\rm in} \ T  \\
  \ u_{ij}(x)  \in {S^{n-p-1}_{-}}&  \ \mbox{\rm if} \ j<_p i  \ \mbox{\rm in} \ T  
  \end{array} \right. \right\} \raisebox{-3mm}{\  .}\]
\end{itemize}
\end{pro}  

\Proof  The equality 
$$\GJ^n_T =
\bigcup_{\tau \hspace{0.4mm}\dom \hspace{0.4mm}T} \FN{\tau} $$
readily follows from the definitions and Theorem \ref{submonadr}. 

Let us prove that  $$X_T = \bigcup_{\tau \hspace{0.4mm}\dom \hspace{0.4mm}T} \FN{\tau}.$$

The inclusion $$\bigcup_{\tau \hspace{0.4mm}\dom \hspace{0.4mm}T} \FN{\tau}\subset X_T $$ is obvious from Lemma \ref{GJtau}.  Let us prove the inverse inclusion.

Let $x\in X_T .$ We have  to prove that $\tau(x)\dom T .$   
Let $i<_p j$ in $\tau(x) .$ Then, by Lemma  \ref{GJtau}, again $u_{ij}(x) \in \stackrel{\scriptscriptstyle o \  \ \ \ \ \ \  }{S^{n-p-1}_{+}} . $ Since $x\in X_T$ two possibilities exist. Either $i<_r j$  in $T$ for some $r$ and then $u_{ij}(x)\in 
{S^{n-r-1}_{+}}$ or $j<_q i$  in $T$ for some $q$ and $u_{ij}(x)\in 
{S^{n-q-1}_{-}} .$

In the  first case we have $x\in \stackrel{\scriptscriptstyle o \  \ \ \ \ \ \  }{S^{n-p-1}_{+}}\cap {S^{n-r-1}_{+}}$ which is possible only if $r\le p .$ In the second case  $x\in \stackrel{\scriptscriptstyle o \  \ \ \ \ \ \  }{S^{n-p-1}_{+}}\cap {S^{n-q-1}_{-}}$ and this is possible only if $q< p.$ So $\tau$ is dominated by $T .$  

\

\Q

\

\begin{corol} The space $\GJ^n_T$ is a closed subspace of $\fm^n_{|T|}$ \end{corol}

Now we  will study the properties of $\GJ^n$ and its cousin $\GJ^n_{\circ} .$

\begin{theorem}\label{Gjcon} The operad $\GJ^n$ is contractible. \end{theorem}

\Proof  For $n=2$ Kontsevich and Soibelman prove contractibility of $X_T$ in \cite{KS}[Proposition 7]. 
Unfortunately the details are left to the reader. So we will follow their idea and prove  it for arbitrary $n$ (but only for $n$-ordinals, not arbitrary complimentary orders as in \cite{KS}). 

As in \cite{KS} we will use  induction on the number of tips of a tree $T .$ If $T$ has only two tips then 
obviously the space $X_T$ is $S^{n-p-1}_{+}$ , where $p$ is the index of  the unique nonempty order on $\{1,2\} .$

Suppose we have proved this theorem for all $T$ with $|T|=[k] .$  Let $T$ be a tree with $|T| = [k+1] , $ and 
let $T'$ be an $n$-tree which is obtained from $T$ by cutting off the most right branch. Then we have an injection 
$$T' \rightarrow T$$
which induces a map $\pi:X_T \rightarrow X_{T'} .$

Let $a\in X_{T'}$ be a point. Let us suppose that  $a$  belongs to the cell $\FN{T'} .$ We will show that the fiber of $\pi$ is homeomorphic to a disk of dimension $n-p .$ Here, $p$ is such that  $k<_p k+1$ in $T .$ 

Since the fiber obviously depends only on configurations of points the labels of which belong to the most right  
 branch of the  tree $T' ,$ then   without loss of generality we can assume that $T'$ 
is a suspension of some $(n-1)$-tree and $k<_{0} k+1$ in $T$. The manifold  $\FN{T'}$ is diffeomorphic to the 
intersection of the space of  configurations of  points 
$$(x^1_1, \ldots,x_1^n),(x^1_2,\ldots,x^n_2),\ldots,(x^1_k,\ldots,x^n_k) \in (\Re^n)^k$$
 which belong to the hyperplane
  $x^1=0 ,$ with $x_1 = 0$ and $x_k^n = 1 $ with the Fox-Neuwirth cell $\fn_{T'} .$ Let $(x_1=0,x_2, \ldots, x_k)$ be the image  of $a$ in this space. In its turn $\FN{T}$ is diffeomorphic  to the 
space of configurations of  points $x_1, \ldots, x_{k+1} $ in $\Re^n$ such that $x_1, \ldots, x_k$ belongs to the 
previous intersection and  $x_{k+1}$ is in the open positive 
halfspace $x^1>0 .$ 

Let $r$ be a sufficiently small number such that the closed balls $B_r(x_i)$ of radius $r$ 
with the centre $x_i, \ 1\le i\le k ,$ do  not intersect   each other. And let $R$ be  sufficiently big 
 such that the union of all $B_r(x_i)$ belongs to the interior of the ball $B_R(0) .$ Let $C(a,r,R)$ be the manifold with corners 
$$C(a,r,R) = \{ x\in \Re^n \ | \  x^1 \ge 0 \ , \  x \in B_R(0)\ , \  x\in \Re^n \setminus (\cup_{i=1}^{k} int(B_r(x_i))) \} .$$

\Example
\vskip10pt \hskip40pt \includegraphics[width=120pt]{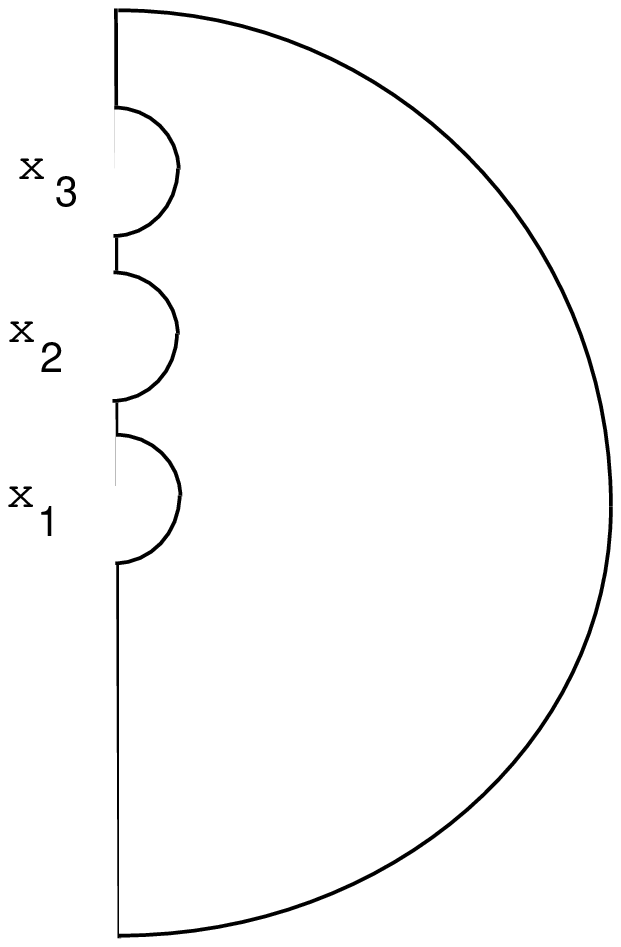}

We claim, that there is a homeomorphism from $C(a,r,R)$ to the fiber of $\pi$ over $a .$    

To construct such a homeomorphism we first choose an $r_1>r$  and $R_1<R$ such that  the $B_{r_1}(x_i)$ still do not
 interstect and their union is still in the interior  of the ball $B_{R_1}(0) .$
For a fixed $1\le i\le k$ let 
$$\psi_i(y) = \frac{r(||y||-r)}{r_1(r_1-r)}(y-x_i)+x_i \ , $$
for $ y\in C(a,r,R) $ such that   $ r \le ||y-x_i|| \le r_1 .$
And let
$$\psi_{\infty}(y) = \frac{R(R_1-R)y}{(R-||y||)R_1} $$
for $ y\in C(a,r,R) $ such that   $ R_1 \le ||y|| < R .$

We can construct a map $F$ from $C(a,r,R)$ to the positive halfspace $x_1\ge 0$ with the following properties:
\begin{itemize}
\item the restriction of $F$ on $C(a,r_1,R_1)$ is a homeomorphism of manifolds with corners;
\item $F(y) = \psi_i(y)$ if $ r \le ||y-x_i|| \le r_1 ;$
\item $F(y) = \psi_{\infty}(y)$ if $ R_1 \le ||y|| < R .$
\end{itemize}

Now we construct a map $F'$ from $C(a,r,R)$ to $\pi^{-1}(a)$ to be equal at a point $y\in C(a,r,R)$ to
$$\lim_{m  \rightarrow  \infty} (u_{ij}(x_1,\ldots,x_k,F(y_m)), d_{ijl}(x_1,\ldots,x_k,F(y_m))$$ 
where $\{y_m\}$ is an arbitrary sequence from the interior of $C(a,r,R)$ which converges to $y.$

Obviously, $F'$ is correctly defined and continuous. It also maps  injectively the  interior  of $C(a,r,R)$ 
 to the intersection of $\pi^{-1}(a)$ with $\FN{T} .$  
  Points from the boundaries of $B(x_i,r)$ are mapped injectively 
 to the points of intersections of $\pi^{-1}(a)$ and $\FN{c}$ 
 where $c$ is the Getzler-Jones cell corresponding to a map of 
trees $\sigma:T\rightarrow T'$, $\sigma(k+1)= i  \  , \ \sigma(j) = j , \ j\le k .$ 
 This intersection is homeomorphic to the hemisphere $ \stackrel{\scriptscriptstyle o \ \ \ }{S^{n-1}_{+}} .$  

Points from the hyperplane $x_1 = 0$ are mapped under $F'$ to the points  of the Getzler-Jones cell  corresponding 
to the configuration $(x_1,\ldots, x_k, y)$ and, finally, points from the outer hemisphere boundary of $C(a,r,R)$
 are mapped to the Getzler-Jones cell corresponding to the 
map of trees $\sigma:T\rightarrow M_0^2 $ , $ \sigma(i) = 1$   if $1\le  i \le k$ and $\sigma(k+1) = 2 .$ 
 
 It is not difficult to check that $F'$ is bijective and so it is a homeomorphism since both spaces are compact. 
 
 Suppose now $a$ belongs to a  generalized  Getzler-Jones cell  in $X_{T'} $  then the fiber $\pi^{-1}(a)$ 
can be glued from the manifolds $C(b_t, r_t,R_t), t\in Vertex(\tau(a)) ,$ by the following inductive procedure.
 First  we construct $C(b_0,r_0,R_0)$ where  $b_0$ is the projection of $a$  to the configuration which decorates 
 the root vertex $v_0$ of $\tau(a) .$ 
 Then we consider  the vertices $v_1,\ldots,v_s$  which can be connected to $v_0$ by exactly one edge. 
We then construct $C(b_l,r_l,R_l) , 1\le l \le s ,$ where $b_l$ is the corresponding projection of $a .$ 
By scaling up or down the configurations $b_l$ if necessary we always can make 
 $R_1 = R_2 = \ldots = R_s = r_0 .$ So we construct the next manifold by gluing the outer hemisphere  $C(b_l,r_l, R_l)$
   to the inner hemisphere of $C(b_0,r_0,R_0)$ in the place $l .$
 
 \
 
 \Example   The following example illustrates the proof. Here 
 $$T = [5]\stackrel{\rho}{\rightarrow} [3] \rightarrow [1] $$  
 $$\rho(1)=\rho(2)=\rho(3)= 1 \ , \rho(4)= 2 \ , \  \rho(5)= 3$$
 and the Getzler-Jones cell corresponds to the map $\sigma: T'\rightarrow S$ of $2$-trees
 $$S=[2]\rightarrow [1]\rightarrow [1] \ , \sigma(1) = \sigma(2) = 1 \ , \  \sigma(3) = \sigma(4) = 2 .$$

 \vskip10pt \hskip40pt \includegraphics[width=220pt]{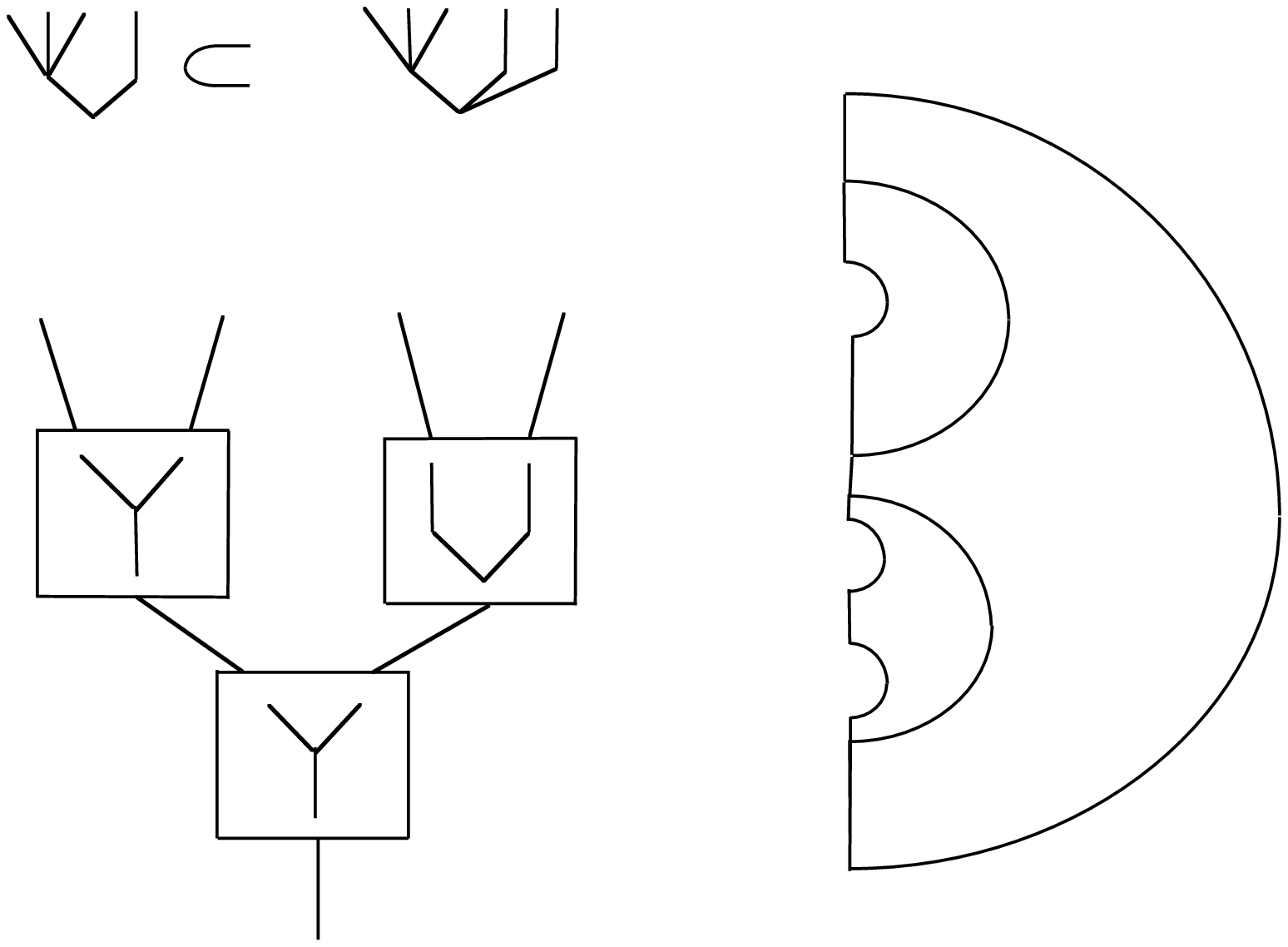}
 
 \
 
 Obviously, this inductive procedure stops after a finite number of steps and produces a contractible manifold.  

\

\Q

\

\begin{theorem}\label{Gjcof} The operad $\GJ_{\circ}^n$ is a cellular object in $\URO(Top) ;$ in particular, it is cofibrant.  \end{theorem}

\Proof  
First recall some terminology  \cite{BergerM}.  We will work in the category of unbased reduced $n$-operads 
and unbased reduced $n$-collections  and we will call them simply $n$-operads and $n$-collections. 
  Let $A\in \URO(Top)$. Let $$V_n: \URO(Top) \rightarrow URColl_n(Top)$$ be the forgetful 
functor with  left adjoint $\URF .$  Let $\omega: V_n(A) \rightarrow  K$ be a cofibration of 
  $n$-collections.  {\it The cellular extension  $A\rightarrow A[\omega]$ generated by $\omega$} 
is the following pushout in $\URO(Top)$

 {\unitlength=1mm
\begin{picture}(60,33)
\put(43,25){\makebox(0,0){\mbox{$\URF V_n (A)$}}}
\put(43,21){\vector(0,-1){10}}
\put(53,26){\shortstack{\mbox{$ $}}}
\put(77,25){\makebox(0,0){\mbox{$\URF(K)$}}}
\put(77,21){\vector(0,-1){10}}
\put(53,25){\vector(1,0){14}}
\put(68,15){\shortstack{\mbox{$$}}}
\put(41,15){\shortstack{\mbox{$ $}}}
\put(43,8){\makebox(0,0){\mbox{$A$}}}
\put(77,8){\makebox(0,0){\mbox{$A[\omega]$}}}
\put(49,8){\vector(1,0){22}}
\put(55,8){\shortstack{\mbox{$ $}}}

\end{picture}}

\noindent where the left vertical map is the counit of the adjunction.

We will prove that $\GJ^n_{\circ}$ can be obtained as a  sequential colimit of   cellular extensions starting from the initial $n$-operad.

Consider a  filtration of the $n$-collection $\FN{\bullet}$ by the following subcollections $\FN{\bullet}(m) \ ,  \ m\ge 0 .$
\[
\FN{T}(m)  =  \left\{ 
\begin{array}{ll}
 \FN{T} & \mbox{if} \ dim(T) \le m   \\
 & \\
\emptyset & \mbox{if} \  dim(T) > m \end{array}
\right.
\]
We have an inclusion  for every $m$
  $$\Phi_m: \URF(\FN{\bullet}(m)) \longrightarrow \URF(\FN{\bullet}) \longrightarrow \GJ^n_{\circ}.$$
The image $\GJ^n_{\circ}(m)$  of this inclusion is closed and, hence, a compact suboperad of $\GJ^n_{\circ} .$ Moreover, 
$$ \GJ^n_{\circ} \simeq \mbox{colim}_m \GJ^n_{\circ}(m) .$$

We want to show that $\GJ^n_{\circ}(m) \subset\GJ^n_{\circ}(m+1)$ is a cellular extension generated by a cofibration.  
Indeed, let us consider the following  $n$-collection $k(m+1).$
\[
  k_{T }(m+1)=  \left\{ 
\begin{array}{ll}
 V_n(\GJ^n_{\circ}(m))_T\cup \FN{T}& \mbox{if} \ dim(T) = m+1   \\
  &  \\
 V_n(\GJ^n_{\circ}(m))_T & \mbox{if} \ dim(T)\ne m+1 \end{array}
\right.
\]

For every $n$-tree $T ,$   $k_{T }(m+1)$ is a closed subspace of $V_n(\GJ^n_{\circ})(m+1)_T .$ Moreover, 
 if $dim(T) = m+1$ then  Proposition \ref{FNclosure} implies that there exists 
a homeomorphism $$S^{ dim(T)} \rightarrow Bd(K_T) = Bd(cl(\FN{T}))$$ such that we have the following  pushout in $Top .$ 

{\unitlength=1mm
\begin{picture}(60,33)
\put(43,24){\makebox(0,0){\mbox{$\coprod\limits_{T,dim(T)=m+1}\hspace{-7mm}S^{ dim(T)}$}}}
\put(40,20){\vector(0,-1){8}}
\put(53,26){\shortstack{\mbox{$ $}}}
\put(79,25){\makebox(0,0){\mbox{$ V_n(\GJ^n_{\circ}(m))_T$}}}
\put(79,21){\vector(0,-1){10}}
\put(56,25){\vector(1,0){9}}
\put(68,15){\shortstack{\mbox{$$}}}
\put(41,15){\shortstack{\mbox{$ $}}}
\put(43,6){\makebox(0,0){\mbox{$\coprod\limits_{T,dim(T)=m+1}\hspace{-7mm}B^{ dim(T)}$}}}
\put(79,7){\makebox(0,0){\mbox{$k_{T }(m+1)$}}}
\put(56,7){\vector(1,0){13}}
\put(55,8){\shortstack{\mbox{$ $}}}
\end{picture}}

\noindent  Hence, the inclusion $$\omega_{m+1}: V_n(\GJ^n_{\circ}(m))\  \subset \  k(m+1)$$ is a cofibration in the category of 
unbased reduced $n$-collections. 

Besides that, we have a map of $n$-operads
$$\URF(k(m+1) \longrightarrow \  \GJ^n_{\circ}(m+1)$$
generated by the obvious map of collections $  k (m+1)\longrightarrow \  \GJ^n_{\circ}(m+1) .$ 

We also can construct another $n$-collection 
\[
  l_{T }(m+1)=  \left\{ 
\begin{array}{ll}
V_n(\URF(\FN{\bullet}(m)))_T\coprod \FN{T}& \mbox{if} \ dim(T) = m+1   \\
& \\
V_n(\URF(\FN{\bullet}(m)))_T & \mbox{if} \ dim(T)\ne m+1 \end{array}
\right.
\]
with a bijective continuous map of $n$-collections 
$$l(m+1) \longrightarrow \ k(m+1)  ,$$ 
 a cofibration 
\begin{equation}\label{coprodL} \varsigma_{m+1}: W_n(\URF(\FN{\bullet}(m)))\  \subset  \   L(m+1) ,\end{equation}
and a map
$$\URF(l(m+1)) \longrightarrow \  \URF(\FN{\bullet}(m+1)) .$$
The cofibration \ref{coprodL} is actually a canonical coprojection into a coproduct. 

All these maps of operads can be organised into a  commutative cube.

{\unitlength=1mm

\begin{picture}(60,55)(-20,-4)

\put(0,25){\makebox(0,0){\mbox{$\URF(V_n(\URF(\FN{\bullet}(m))))$}}}
\put(10,22){\vector(0,-1){14}}
\put(12,15){\shortstack{\mbox{$ $}}}

\put(24,25){\vector(1,0){15}}

\put(23,26){\shortstack{\mbox{$ $}}}

\put(53,25){\makebox(0,0){\mbox{$\URF(l(m+1))$}}}
\put(51,22){\vector(0,-1){14}}

\put(57,21){\shortstack{\mbox{$ $}}}

\put(57,28){\shortstack{\mbox{\small $ $}}}

\put(5,5){\makebox(0,0){\mbox{$\URF(\FN{\bullet}(m))$}}}

\put(19,5){\vector(1,0){14}}

\put(23,6){\shortstack{\mbox{$ $}}}

\put(51,5){\makebox(0,0){\mbox{$ \URF(\FN{\bullet}(m+1))$}}}

\put(57,1){\shortstack{\mbox{$ $}}}

\put(57,8){\shortstack{\mbox{\small $ $}}}

\put(-3,3){\begin{picture}(50,30)


\put(25,35){\makebox(0,0){\mbox{$\URF(V_n(\GJ^n_{\circ}(m)))$}}}

\put(30,32){\line(0,-1){9.5}}
\put(30,21.2){\vector(0,-1){3}}

\put(32,25){\shortstack{\mbox{$ $}}}

\put(42,35){\vector(1,0){15}}

\put(43,36){\shortstack{\mbox{$ $}}}

\put(71,35){\makebox(0,0){\mbox{$\URF(l(m+1))$}}}
\put(72,32){\vector(0,-1){14}}

\put(77,31){\shortstack{\mbox{$ $}}}

\put(77,38){\shortstack{\mbox{\small $ $}}}

\put(30,15){\makebox(0,0){\mbox{$ \GJ^n_{\circ}(m)$}}}

\put(55,15){\vector(1,0){4.4}}
\put(40,15){\line(1,0){12.4}}

\put(43,16){\shortstack{\mbox{$ $}}}

\put(72,15){\makebox(0,0){\mbox{$  \GJ^n_{\circ}(m+1)$}}}

\put(77,11){\shortstack{\mbox{$ $}}}

\put(77,18){\shortstack{\mbox{\small $ $}}}

\end{picture}}

\put(13,28){\vector(1,1){7}}

\put(54,28){\vector(1,1){7}}

\put(13,8){\vector(1,1){7}}

\put(54,8){\vector(1,1){7}}

\end{picture}}

\noindent The front square of this cube is  a pushout in the category of $n$-operads because $\varsigma_{m+1}$ 
is a coprojection and the free operad functor preserves colimits.   All   maps from the front square to the back square are continuous bijections.
 So the induced morphism from  $ \URF(\FN{\bullet}(m+1))$ to the pushout $P$ of the back square must be a continuous bijection as well,
 and, hence, the induced map  $\alpha: P\rightarrow  \GJ^n_{\circ}(m+1)$ is a continuous bijection.  

Moreover, the vertical map $\URF(l(m+1)) \rightarrow \URF(\FN{\bullet}(m+1))$ admits a section. It follows that 
$\URF(k(m+1))_T\rightarrow P_T$ is epi. Since $\URF(l(m+1))_T$ is a finite coproduct of compact spaces the space $P_T$ is compact. 
Therefore,  $\alpha$ is an isomorphism and we have proved our theorem.

\

\Q

\

By definition we have an inclusion $\GJ^n \rightarrow Des_n(\fm^n),$ and, by adjunction, a canonical  morphism
$$Sym_n(\GJ^n)\rightarrow \fm^n .$$  

\begin{theorem}\label{SymGJ} This canonical  map from $Sym_n(\GJ^n)$ to $\fm^n$ is an isomorphism. The analogous result holds in the unbased case. 
\end{theorem}  

\Proof  
The canonical map 
$$ \RF_n(\FN{\bullet})\rightarrow \GJ^n $$
is a continuous bijection. Hence, after application of $Sym_n$ to this map we have a  composite of continuous bijections
$$\RSF(C_n(\FN{\bullet})) \rightarrow  Sym_n(\RF_n(\FN{\bullet}))\rightarrow Sym_n(\GJ) .$$

We also have  the following commutative diagram

 {\unitlength=1mm
\begin{picture}(60,33)
\put(15,25){\makebox(0,0){\mbox{$\RSF(C_n(\FN{\bullet}))_k$}}}
\put(15,21){\vector(0,-1){10}}
\put(55,26){\shortstack{\mbox{$ $}}}
\put(80,25){\makebox(0,0){\mbox{$Sym_n(\RF_n(\FN{\bullet})) $}}}
\put(80,21){\vector(0,-1){10}}
\put(37,25){\vector(1,0){20}}
\put(68,15){\shortstack{\mbox{$ $}}}
\put(42,15){\shortstack{\mbox{$ $}}}
\put(15,8){\makebox(0,0){\mbox{$\F_s(\mod{k})$}}}
\put(80,8){\makebox(0,0){\mbox{$Sym_n(\GJ^n) $}}}
\put(30,8){\vector(1,0){10}}
\put(48,8){\makebox(0,0){\mbox{$\fm^n$}}}
\put(65,8){\vector(-1,0){10}}
\put(55,8){\shortstack{\mbox{$ $}}}
\end{picture}}

\noindent where the composite of the left vertical and horizontal maps is a continuous bijection by \ref{GJdecomp}.
Hence, the canonical map from $Sym_n(\GJ^n)$ to $\fm^n$ is a continuous bijection, as well.
By formula \ref{formular} the space  $Sym_n(\GJ^n)_k$ is a finite colimit of compact spaces, hence,
 it is compact and therefore the bijection between $Sym_n(\GJ^n)$ and $\fm^n$ is a homeomorphism.  

\

\Q

\section{$n$-operads and $E_n$-operads}

We finally are able to consider some applications of the results obtained.

\begin{theorem}\label{fmcof} The operad $\fm^n_{\circ}$ is a cellular object in the category of unbased reduced symmetric operads. In particular, it is  cofibrant. \end{theorem}
\Proof This is an easy consequence of the fact that $Sym_n$ is a left Quillen functor. 

\

\Q

\

\Remark The cofibrantness of $\fm^n_{\circ}$ was first claimed in \cite{GJ} without a proof. To the best of our knowledge the proof first appeared in \cite{Salvatore} and uses a comparison between $\fm^n_{\circ}$ and its Boardman-Vogt $W$-construction. 

\

For a categorical operad $A$ we will denote by $|A|$ the geometric realisation of the nerve of $A .$

\begin{theorem}
\label{rhnfmnstrong} The operad $|\urh|$ is strongly homotopy equivalent to
$\fm^n_{\circ}$ . \end{theorem}

\Proof Both $n$-operads $|\URH|$ and $\GJ^n_{\circ}$ are fibrant, cofibrant and contractible. 
Hence, they are strongly homotopy equivalent. So are their images under $Sym_n .$ 

\

\Q

\

\Remark The Getzler-Jones decomposition of $\fm^n_{\circ}$ is not a regular $CW$-decomposition so we cannot
 take the poset of its cells and form a categorical operad as was proposed in \cite{GJ}. 
However, the previous corollary shows that $\urh$ is an appropriate  substitute for this nonexistent poset operad.   

\

\begin{theorem}\label{rhnfmn} There is a chain of weak operadic equivalences between $|\rh^n|$ and $\fm^n .$ \end{theorem}

\Proof The operad $|\RH^n|$ is the bar construction on the terminal reduced $n$-operad. So, we have a zig-zag
$$B(\RF_n,\RF_n,1)\longleftarrow B(\RF_n,\RF_n,\GJ^n)\longrightarrow
\GJ^n $$
which, after symmetrization,  gives the following zig-zag of morphisms of symmetric operads
 \begin{equation}\label{zigzag}|\rh^n|\longleftarrow Sym_n B(\RF_n,\RF_n,\GJ^n)\longrightarrow
\fm^n \end{equation}
by Theorem \ref{Symnerve}.
If we forget about nullary operations, we get the  zig-zag of weak operadic equivalences
$$|\urh|\longleftarrow Sym_n B(\URF,\URF,\GJ^n_{\circ})\longrightarrow
\fm^n_{\circ}  $$
 since it can be obtained by symmerisation from the zig-zag of equivalences of fibrant cofibrant operads
$$|\URH|\longleftarrow B(\URF,\URF,\GJ^n_{\circ})\longrightarrow
\GJ^n_{\circ} . $$
Hence, the zig-zag \ref{zigzag} consists of weak equivalences.  

\

\Q

\

 Recall that the iterated monoidal category operad $\M^n$ contains an internal $n$-operad \cite{BEH}.
It is  easy to see that $\M^n$ is a reduced categorical symmetric operad and its internal $n$-operad is also a reduced internal $n$-operad.
So we have a canonical map of operads
$$k^n: \rh^n\rightarrow \M^n .$$
  Recall also that the inclusion of the Milgram poset $$(\J^n_k)^{op}\rightarrow \M^n_k$$ induces a weak equivalence on the nerves \cite{BFSV}. 

The following theorem provides an alternative proof of the Theorem  \ref{rhnfmn}.

\begin{theorem}\label{Nk} The  map of operads $N(k^n): N(\rh^n) \rightarrow N(\M^n)$ is a weak equivalence. \end{theorem}

\Proof   
By Theorems \ref{Symnerve} and \ref{formular},  we have 
$$N(\rh^n_k) \simeq   \mbox{\rm co}\hspace{-2.4mm}\!\lim\limits_{\mbox{\rm (}\J^{n}_k\mbox{\rm )}^{\raisebox{0.1mm}{\scriptsize op}}}
 \widetilde{{N(\RH^n)}_T}                                                                     .$$
It is trivial to check that $\widetilde{{N(RH^n)}_{(-)}}$ is a Reedy cofibrant functor on $(\J^n_k)^{op} $ because  
$$ \mbox{\rm co}\hspace{-5.3mm}\!\lim\limits_{T\leftarrow T' , T\ne T' } \widetilde{{N(\RH^n)}_{T'}}  \longrightarrow  \widetilde{{N(\RH^n)}_T} = N(\RH^n)_T$$
is a monomorphism. 
Therefore, 
$$\mbox{\rm hoco}\hspace{-2.4mm}\!\lim\limits_{\mbox{\rm (}\J^{n}_k\mbox{\rm )}^{\raisebox{0.1mm}{\scriptsize op}}} \widetilde{{N(\RH^n)}_T}
 \longrightarrow \mbox{\rm co}\hspace{-2.4mm}\!\lim\limits_{\mbox{\rm (}\J^{n}_k\mbox{\rm )}^{\raisebox{0.1mm}{\scriptsize op}}}
 \widetilde{{N(\RH^n)}_T}                                                                     $$
is a weak equivalence.   But $N(\RH^n_T)$ is  contractible, hence, we have a weak equivalence
$$\mbox{\rm hoco}\hspace{-2.4mm}\!\lim\limits_{\mbox{\rm (}\J^{n}_k\mbox{\rm )}^{\raisebox{0.1mm}{\scriptsize op}}} {N(\J^n_k / T)}  \longrightarrow  \mbox{\rm hoco}\hspace{-2.4mm}\!\lim\limits_{\mbox{\rm (}\J^{n}_k\mbox{\rm )}^{\raisebox{0.1mm}{\scriptsize op}}}\widetilde{{N(\RH^n)}_T} .$$  
So in   the commutative square  
{\unitlength=1mm

\begin{picture}(60,33)
\put(35,23.5){\makebox(0,0){\mbox{$ \mbox{\rm hoco}\hspace{-2.4mm}\!\lim\limits_{\mbox{\rm (}\J^{n}_k\mbox{\rm )}^{\raisebox{0.1mm}{\scriptsize op}}} {N(\J^n_k / T)}$}}}
\put(37,21){\vector(0,-1){10}}
\put(53,26){\shortstack{\mbox{$ $}}}

\put(72,25){\makebox(0,0){\mbox{$N(\J^n_k)$}}}
\put(72,21){\vector(0,-1){10}}

\put(52,25){\vector(1,0){12}}

\put(35,7.5){\makebox(0,0){\mbox{$\mbox{\rm hoco}\hspace{-2.4mm}\!\lim\limits_{\mbox{\rm (}\J^{n}_k\mbox{\rm )}^{\raisebox{0.1mm}{\scriptsize op}}} \widetilde{{N(\RH^n)}_T}$}}}

\put(72,8.5){\makebox(0,0){\mbox{$N(\rh^n_k)$}}}

\put(52,8.5){\vector(1,0){12}}

\end{picture}}
                                         
\noindent both horizontal and left vertical arrows are weak equivalences and, therefore,  the right vertical arrow  
is a weak equivalence. 

Finally, the functor $(\J^n_k)^{op} \rightarrow \M^n_k$   is factorised as 
$$(\J^n_k)^{op} \longrightarrow \rh^n_k \longrightarrow \M^n_k$$
and the statement of our theorem follows.

\

\Q

\

\begin{cor} The canonical map $|\urh| \rightarrow |\M^n_{\circ}|$ provides  a cofibrant replacement 
for $|\M^n_{\circ}|$ where $\M^n_{\circ}$ is the unbased version of $\M^n .$ \end{cor} 

There is a canonical map  $\psi: \ph^n\rightarrow \rh^n$ since $\rh^n$ contains an internal pruned operad.
  We actually can restrict $\psi$  to a subcategory of $\ph^n$ which consists of objects which do not contain degenerate decorations. Obviously, the last subcategory is a deformation retract of $\ph^n .$ The deformation retraction is given by the following composite (dropping off of dead leaves):
$$ \mu(T; z^nU_0, x_1, \ldots, x_k) \rightarrow  \mu(T; z^nU_0, \mu(U_n, x_1), \ldots, \mu_(U_n, x_k))= $$
$$= 
\mu(\mu(T; z^nU_0,  U_n, \ldots, U_n); x_1, \ldots, x_k) \rightarrow \mu(T'; z^nU_0, x_1, \ldots, x_k) ,$$
where $T'$ is obtained from $T$ by dropping off a branch which corresponds to the degenerate fiber.

So we will consider this subcategory as the domain of $\psi$ but abusing notation will call it $\ph^n .$ 

Let  $\tau\in \rh^n_k $ be a labelled tree decorated by pruned $n$-trees.   Obviously, an object from the fiber of $\psi$ over $\tau$ will be  a labelled planar tree  decorated by pruned $n$-trees
such that the reduction of it (i.e. the deletion   of all   vertices  of valency $2$ together with their decorations)  is  
$\tau .$ The morphisms are insertion of a vertex of valency $2$    decorated by $U_n$ and deletion of such decorations. 
The following lemma is obvious from this description. 

\begin{lemma}\label{fiber}  The reduced decorated tree $\tau$ considered as an object of the fiber $\psi^{-1}(\tau)$ is a terminal object of this fiber. \end{lemma}

\begin{lemma}\label{section} The operadic morphism $\psi$ has a (nonoperadic) section 
$$s: \rh^n \rightarrow \ph^n$$
which maps an object $\tau$ to the terminal object $\tau$ in $\psi^{-1} .$ 
\end{lemma}

\Proof To define $s$ on morphisms we have to define it on generators and then check correctness.
We define it on a generator corresponding to a surjection $\sigma:T  \rightarrow S$
equal to the composite 
$$\mu(S; T^{(p)}_1, \ldots, T^{(p)}_k) \rightarrow \mu(S; T_1^{(p)}, \ldots, U_n, \ldots, T_k^{(p)}) \rightarrow T$$
where we insert $e\rightarrow U_n$ in the places where $\sigma$   has one-tip fibers. It is 
routine to check the relations. 

\

\Q

\

From these two lemmas we have 

\begin{theorem}\label{phnrhnmn} The operadic functor $\psi$ induces a weak equivalence of simplicial operads
$$N(\psi): N(\ph^n) \rightarrow N(\rh^n);$$
so  all three operads $N(\ph^n), N(\rh^n), N(\M^n)$ are $E_n$-operads.
\end{theorem}

A  pruned  topological $n$-operad  $A$  will be called  {\it contractible}  provided the unique map to the terminal $n$-operad is a weak equivalence  i.e.  every $A_T$ is a contractible topological space.
  
\begin{theorem}\label{Enspace}
Let $A$ be a contractible pruned $n$-operad in the category of compactly generated Hausdorff spaces  such that every $A_T$ 
is a cofibrant topological space   and let $X$ be an algebra
of $A$. Then $X$ has a structure of an $E_n$-space, so up to group completion $X$ is an $n$-fold loop space.
\end{theorem}

\Proof   Since $X$ is an algebra of $A$ it is also an algebra  of $Sym_n (A)$ and, 
by  Theorem \ref{formulap},  $Sym_n(A) \simeq
\
\colim \ \tilde{A}_k \ . $ It is not hard to check that the sequence $ \hocolim \
\tilde{A}_{\ast} \ $ has the structure of an operad and, moreover, the canonical map 
\begin{equation}\label{EHF} \ \ \ \ \ \ \ \
\ \ \ \ \ \ \ \ \ \
\ \ \ \ \ \ \ \ \ \ \ \ \ \
\hocolim \ \tilde{A}_{k} \ \rightarrow \ \
\colim \ \tilde{A}_{k} \ \ \ \
\ \ \ \ \ \ \ \ \ \ \
\ \ \ \ \ \ \ \ \ \ \ \ \ \ \ \ \ 
\end{equation}
 is operadic.  But $\hocolim \
\tilde{A}_k
\
$ has the same homotopy type as $\ph^n_k$ because of contractibility of $A_T$. So $X$ is an algebra of a $E_n$-operad
$\hocolim \
\tilde{A}_k .$

\

\Q

\

\Example  It is still possible for $X$ from the previous theorem to be an 
$E_m$-space for $m>n$. For example,
 if $B$ is any
$E_{\infty}$-operad, then $A=Des_n (B)$ is contractible but $Sym_n (A) \simeq B $.

\

A similar theorem holds for reduced $n$-operads in the category $Ch(R)$ of  chain complexes over a commutative ring with unit $R .$ 

Let $M$ be a reduced  $n$-operad in $Ch(R)$ which has $M_T=R$ for every pruned tree. 
A  reduced $n$-operad in $Ch(R)$ equipped with an augmentation $A\rightarrow M $  will be called  {\it contractible}  provided its augmentation is a weak equivalence.  

The method used in the proof of   Theorem \ref{Enspace} can be used without change to prove the following. 
\begin{theorem}\label{Enchain}
Let $A$ be a contractible reduced  $n$-operad in $Ch(R)$  such that $A_T$ is a chain complex of projective $R$-modules for every $T .$ Let $X$ be an algebra
of $A$. Then $X$ admits an action of a symmetric reduced  operad  weakly equivalent to the operad of $R$-chains of the little disc operad.  \end{theorem}

\  

Finally,  Theorems  \ref{Gjcof} and \ref{SymGJ} imply that we can obtain  the full solution of the  coherence problem of $n$-fold 
loop spaces in the spirit of Stasheff's original work \cite{St1,St2} using cells $K_T$ for higher $n$-trees instead of associahedra.
 This was first claimed by Getzler and Jones but 
some doubts appeared  since Tamarkin came up with his counterexample.  Our Theorems \ref{Gjcof}  and \ref{SymGJ} do  confirm 
that the $\fm^n_{\circ}$-algebra structure  is equivalent to the existence of a sequence of inductive 
  extensions of  higher homotopies  from the boundaries of $K_T$ to their interior  even though  $\GJ_T$ is not always a PL-ball.   The exact  combinatorics of $K_T$ will be discussed in  a future   paper. 
Here we just give a few examples of the manifolds $K_T$ which can be drawn on  paper.

If $n=1$ and $T= [m]$ then $K_T = K_m ,$ the Stasheff associahedron as  was said before. 
 
 If $n=2$   we have two one-dimensional manifolds commonly known as associator and braiding.

 \vskip20pt \hskip60pt \includegraphics[width=180pt]{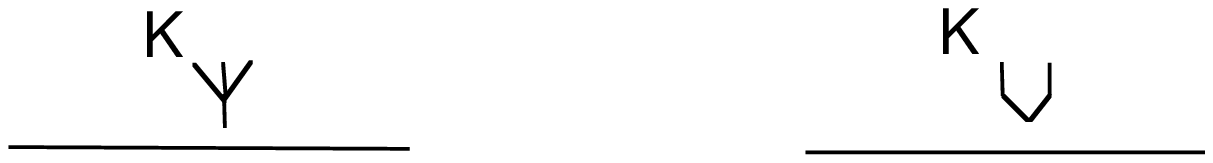}

\

We also have three $2$-dimensional polytopes: a pentagon and two hexagons,   and  two $3$-dimensional polytopes, the associahedron $K_5$
 and a $3$-dimensional polytope
well known   in the theory of Yang-Baxter operators.  
\

 \vskip20pt \hskip50pt \includegraphics[width=250pt]{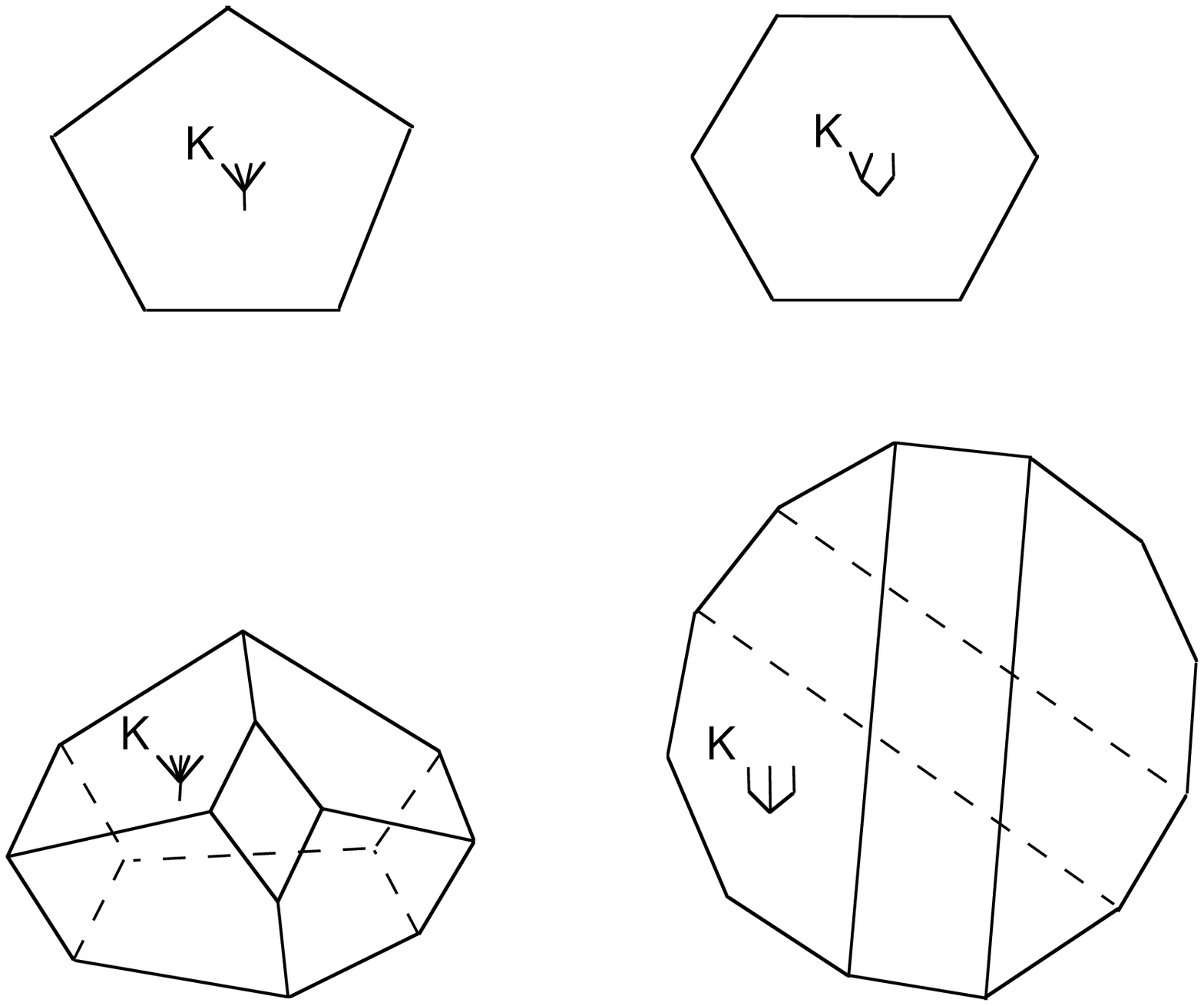}
 
\
  
  Finally we have three $3$-dimensional polytopes (there should be a copy of the first one we do not show here)  which 
 can be found  in the paper of Bar-Natan \cite{Bar-Natan}.  

\

\vskip20pt \hskip50pt \includegraphics[width=170pt]{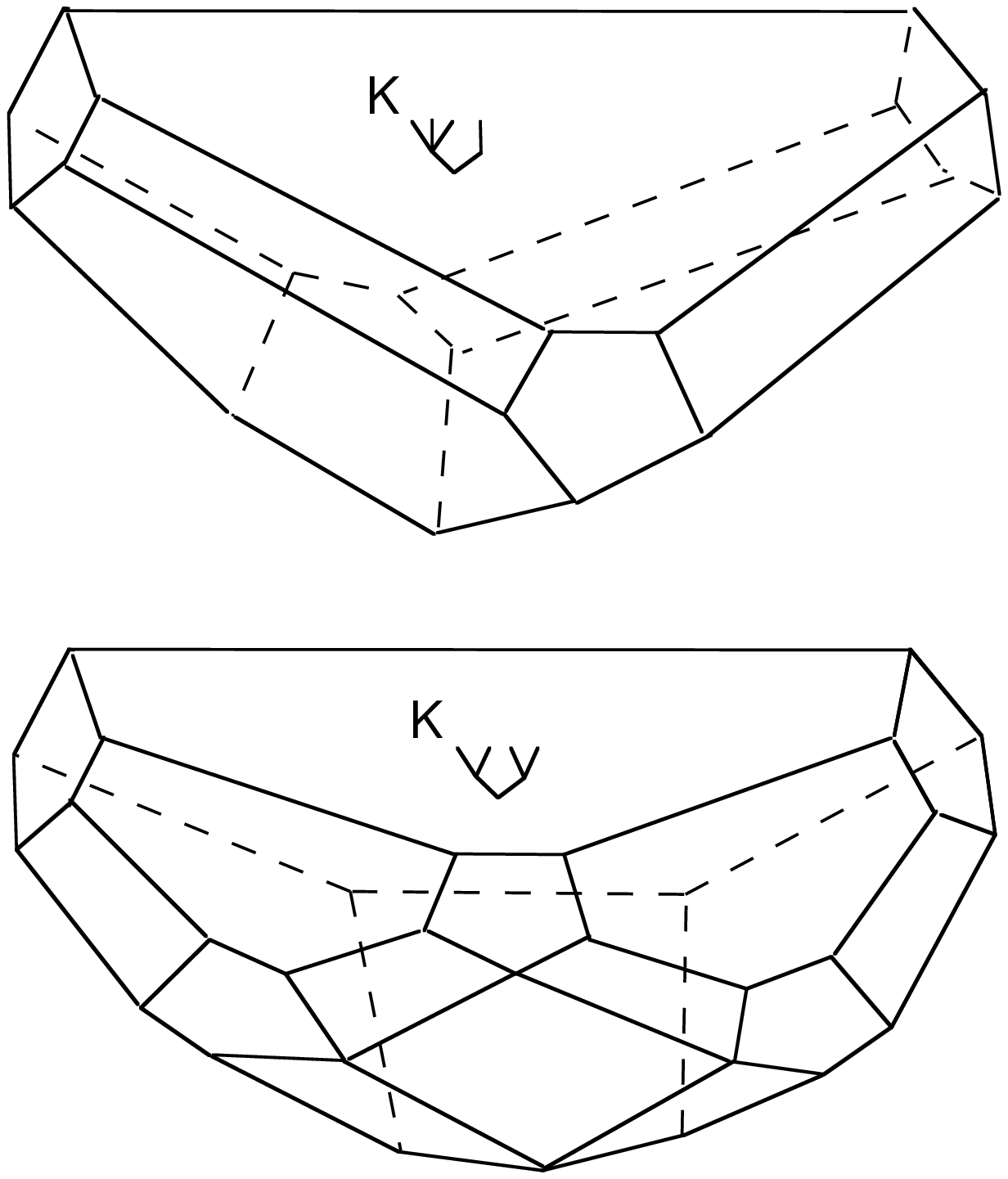}

\

It is recommended to the reader to try to draw up all the polytopes $K_T$ of dimension less than or equal to $3$ with $n=3$ (there will be only three  of them which are different from the above polytopes) and $n=4$ (only one new polytope).

\section{Swiss-Cheese type $n$-operads and their symmetrisation}

This is a short section in which we define reduced (both based and unbased) Swiss-Cheese type $n$-operads and show that the techniques 
developed in the previous sections are powerful enough to easily carry out the main results for classical operads to the case of Swiss-Cheese type operads.

The Swiss-Cheese type operads ({\it SC-operads or SC type operads} for short) were introduced by Voronov in \cite{VoronovSC} with   motivation to describe a finite dimensional model of the moduli space of genus-zero Riemann surfaces from open-closed string theory. 
The importance of this class of operads was also understood  by Kontsevich in \cite{K} who explained that the category of algebras of  Voronov's Swiss-Cheese operad is a natural place for developing a theory of deformation  complexes and its  higher dimensional generalisation. 

Recall first \cite{VoronovSC} that a symmetric reduced SC-operad in a symmetric monoidal category $V$ is a special symmetric coloured operad with two colours and consists      
of a collection $A_{k,l}$ of objects of $V$ indexed by pairs of natural numbers with an action of the product of the symmetric groups $\Sigma_k\times\Sigma_l .$ We also require that $A_{0,0} = A_{1,0}=A_{0,1} = I$ 
for a reduced operad. The operadic multiplication is represented by a family of morphisms
$$A_{\textstyle k,l}\otimes (A_{\textstyle  a_1, b_1}\otimes \ldots \otimes A_{\textstyle a_k,b_k}) \otimes (A_{\textstyle 0, c_1 }\otimes \ldots A_{\textstyle 0,c_l})  \longrightarrow$$
 $$ \longrightarrow A_{\textstyle a_1+\ldots+ a_k \  ,  \  b_1+\ldots+ b_k + c_1+\ldots+ c_l}$$
which must satisfy some natural associativity, equivariancy and unitarity conditions.

For the unbased case (which is actually considered in \cite{VoronovSC, K}) we forget about the space
$A_{0,0} .$ 

Notice that we have two symmetric operads here: $A_{\bullet,0}$ and $A_{0,\bullet} .$  An algebra of such an operad consists 
of a pair of objects $(X_1,X_2)$ such that $A_{\bullet,0}$ acts on $X_1$ and $A_{0,\bullet} $ acts on $X_2 $ and all the other spaces $A_{k,l}$ provide the interplay between these two actions. 

The main example of an unbased reduced symmetric SC-operad  is Voronov's Swiss Cheese operad $\mbox{\bf SC}^n$ whose $(k,l)$-space is the space of $l$ disjoint $n$-disks and $k$ disjoint $n$-semidisks in a big $n$-semidisk 
 $$B^n_{+} = \{(x_1,\ldots,x_n) \in \Re^n \ | \  x_1^2 + \ldots + x_n^2 \le 1 \ , \ x_1 \ge 0 \}$$   
 Notice that our presentation of $\mbox{\bf SC}^n_{k,l}$ is different from Voronov's picture in \cite{VoronovSC}. He defines 
 $B^n_{+}$ by requiring  $x_n \ge 0 .$  It is more convenient for us to ask $x_1\ge 0 ,$ however,  because it makes 
the relations between $\mbox{\bf SC}_{k,l}$ and the Getzler-Jones  decomposition more evident. 
 
\begin{defin} A coloured pruned $n$-tree is a pruned   $n$-tree $T$ equipped with a map of trees 
$$\xi_T: T\longrightarrow M_0^2$$
such that the induced $n$-ordinal structure on  $\xi^{-1}(1)$  is a suspension over $(n-1)$-tree or a degenerate $n$-tree. \end{defin}

One can think of a coloured $n$-tree as an $n$-tree with a distinguished (coloured by $1$) most left branch, which can be empty  however. 

\begin{defin} A coloured morphism between  coloured trees $\sigma:T \rightarrow S$   is a morphism of underlying    trees such that  $\xi_S(\sigma(a)) = 1$ for $a\in  \xi^{-1}(1) .$ \end{defin}
In other words   $\sigma  $ sends the distinguished branch of $T$  to the distinguished branch of $S .$
Every coloured morphism between $n$-trees can be restricted to the distinguished branches of these trees.
 Since these branches are suspensions of $n$-trees, one can consider the restriction of a morphism as a morphism between $(n-1)$-trees.  

A  coloured morphism between coloured $n$-trees is {\it an injection} if the underlying morphism of trees is an injection.
A  coloured morphism between coloured $n$-trees is {\it a surjection} if the underlying morphism of trees is a surjection.
One can take a fiber of  a coloured surjection between $n$-trees and then the   fibers  obtain  canonical colouring.

\begin{defin} A reduced SC type $n$-operad in a symmetric monoidal category $V$ is a  collection  $A_T\in V $ where $T$ runs over all coloured pruned $n$-trees such that $A_{U^{1}_n} = A_{U^{2}_n}=I$ and $A_{z^nU_0}= I,$ where $U^1_n$ is a linear tree with its unique tip coloured by  $1$ and $U^2_n$  is a linear tree with its tip coloured by  $2 .$ This collection is equipped  with 
 a morphism
$$m_{\sigma}: A_S\otimes A_{T_1^{(p)}}\otimes \ldots \otimes A_{T_k^{(p)}} \rightarrow A_T $$
for every coloured morphism of trees $\sigma: T\rightarrow S$ between coloured pruned $n$-trees.   
They must satisfy the obvious associativity and unitarity conditions. \

For an unbased reduced SC-operad we use collections without degenerate trees and define multiplication only with respect to surjections of coloured trees. 

\end{defin}

Immediately from the definition we see that as in the symmetric case an SC type $n$-operad $A$  gives rise to 
two operads $A_{z^n U_n,\bullet}$ and $A_{\bullet , z^nU_n}$ 
 if we consider its restriction to the $n$-trees with no  branches with colour $1$  
 or to the trees with  no branches coloured $2 .$ But unlike in the symmetric case these two operads are of different types: $A_{z^n U_n,\bullet}$  is an $n$-operad but $A_{\bullet , z^nU_n}$  is  an $(n-1)$-operad. 

The reduced SC type $n$-operads form a category and we have a desymmetrisation functor from 
symmetric reduced SC-operads to the reduced SC type  $n$-operads. All the machinery of 
Section \ref{reducedoperads} is then applicable to the Swiss-Cheese case. We denote by 
$SCSym_n$ the corresponding functor of symmetrisation and by $\scrh^n$ the free reduced symmetric SC-operad freely generated by an internal reduced SC type  $n$-operad. 
We have the following analogue of  Theorem \ref{formular}

\begin{theorem}\label{formularsc} Let $A$ be a cocomplete reduced symmmetric   categorical $n$-operad of SC type  and
$x$ be an internal reduced SC type
$n$-operad in $A$ then 
$$
(SCSym_n(a))_{k,l} \ \simeq \ \mbox{\rm co}\hspace{-2.8mm}\!\lim\limits_{\scrh^{n}_{k,l}} \tilde{x}_{k,l}
                                                   $$
where $\tilde{x}_{k,l}:\scrh^{n}_{k,l} \rightarrow A_{k,l}$ is the operadic functor representing the  operad $x .$

In addition,
$$(SCSym_n(a))_{0,\bullet} \simeq Sym_n(A_{z^n U_n,\bullet}) \  ,$$
$$(SCSym_n(a))_{\bullet,0} \simeq  Sym_{n-1} (A_{\bullet , z^n U_n}) .$$

The analogous formula holds in the  unbased case. \end{theorem} 

The unbased Swiss-Cheese analogue of the Fulton-Macpherson operad $\scfm^n_{\circ}$ was defined by Voronov in \cite{VoronovSC}. 
The definition of its  based version $\scfm^n$ is obvious. It is also obvious now how to define the analogue of
the  Getzler and Jones operad $\SCGJ^n .$  A little gift for us is that, for a pruned coloured $n$-tree $T ,$ 
the spaces $\SCGJ_T^n$ and $\GJ_T^n$ coincide. So we have the proof of our next theorem almost  for free:  
\begin{theorem} \label{SymSCGJ}The following analogues of  Theorems \ref{Gjcon},\ref{Gjcof},\ref{SymGJ}, \ref{fmcof}, 
\ref{rhnfmnstrong}, \ref{rhnfmn}   hold: \begin{itemize}  
\item The operad $\SCGJ^n$ is contractible and $\SCGJ^n_{\circ}$ is a cellular and contractible 
unbased reduced SC type $n$-operad.
\item The canonical  map  
$$SCSym_n(\SCGJ^n) \longrightarrow \scfm^n$$ is an isomorphism. The analogous result holds in the  unbased case. 
\item The operad $\scfm^n_{\circ}$ is cellular.
\item The operad $|\scrh^n_{\circ}|$ is strongly homotopy equivalent to
$\scfm^n_{\circ} .$ 
\item  There is a chain of weak operadic equivalences between $|\scrh^n|$ and $\scfm^n .$
\end{itemize}
 \end{theorem}

 It is not hard now to give the definitions of the endomorphism SC $n$-operad  of a pair of objects $X_1,X_2\in V$ and of the 
 category of algebras of an SC $n$-operad.  
 The desymmetrisation  functor preserves the endomorphism operad and, hence, symmetrisation preserves the category of algebras. 
 Analogously to   Theorems  \ref{Enspace} and \ref{Enchain} we have
  
  \begin{theorem}\label{EnspaceSC}
Let $A$ be a contractible reduced SC type  $n$-operad in the category of compactly generated Hausdorff spaces  such that every $A_T$ is a cofibrant topological space and let $(X_1,X_2)$ be an algebra
of $A$. Then 
 $(X_1,X_2)$ admits an action of a symmetric  SC-operad weakly equivalent to  Voronov's reduced  Swiss-Cheese operad. 
\end{theorem}
  
  \begin{theorem}\label{EnchainSC}
Let $A$ be a contractible reduced SC type  $n$-operad in $Ch(R)$  such that $A_T$ is a chain complex of projective $R$-modules for every $T .$
 Let $(X_1,X_2)$ be an algebra
of $A$. Then $(X_1,X_2)$ admits an action of a symmetric reduced  SC-operad  weakly equivalent to the operad of $R$-chains of the Voronov's reduced  Swiss-Cheese operad.  \end{theorem}
  
 Analogous theorems hold in the unbased case.

\end{document}